\documentclass[pdflatex,sn-mathphys-num]{sn-jnl}

\usepackage{graphicx}%
\usepackage{multirow}%
\usepackage{amsmath,amssymb,amsfonts}%
\usepackage{amsthm}%
\usepackage{mathrsfs}%
\usepackage[title]{appendix}%
\usepackage{xcolor}%
\usepackage{textcomp}%
\usepackage{manyfoot}%
\usepackage{booktabs}%
\usepackage{algorithm}%
\usepackage{algorithmicx}%
\usepackage{algpseudocode}%
\usepackage{listings}%
\theoremstyle{thmstyleone}%
\newtheorem{theorem}{Theorem}
\newtheorem{proposition}[theorem]{Proposition}%

\theoremstyle{thmstyletwo}%
\newtheorem{remark}{Remark}%

\theoremstyle{thmstylethree}%
\newtheorem{definition}{Definition}%

\newcommand{\wL}{\overset{\leftarrow}{w}}
\newcommand{\wR}{\overset{\rightarrow}{w}}

\newcommand{\wLx}{\overset{\leftarrow}{w}_{x}}
\newcommand{\wRx}{\overset{\rightarrow}{w}_{x}}

\newcommand{\wLy}{\overset{\leftarrow}{w}_{y}}
\newcommand{\wRy}{\overset{\rightarrow}{w}_{y}}

\newcommand{\qu}{q_u}
\newcommand{\qv}{q_v}

\usepackage{bm}

\begin{document}

\title[A SEMI-IMPLICIT EWB RELAXATION SCHEME FOR 2D SW]{A SEMI-IMPLICIT EXACTLY WELL-BALANCED RELAXATION SCHEME FOR THE 2D SHALLOW WATER SYSTEM}

\author*[1]{\fnm{JOS\'E IGNACIO} \sur{RAM\'IREZ-FUENTES}}\email{nramirez@uma.es}

\author[1]{\fnm{CIPRIANO} \sur{ESCALANTE}}\email{escalante@uma.es}

\author[1]{\fnm{TOM\'AS} \sur{MORALES DE LUNA}}\email{tmorales@uma.es}

\author[1]{\fnm{MAR\'IA LUZ} \sur{MU\~{N}OZ}-RUIZ}\email{mlmunoz@uma.es}

\affil*[1]{\orgdiv{Department of Mathematical Analysis, Statistics and Operations Research, and Applied Mathematics}, \orgname{University of Malaga}, \orgaddress{\street{Campus de Teatinos s/n}, \city{Malaga}, \postcode{29071}, \country{Spain}}}

\abstract{We propose a semi-implicit, exactly well-balanced finite volume scheme for the two-dimensional shallow water equations over a non-flat bottom, aimed at the low-Froude-number regime. The method combines a Suliciu-type relaxation of the pressure with a flux splitting that separates the fast acoustic waves from the slow transport ones: the pressure subsystem is integrated implicitly and the transport subsystem explicitly, so that the time step is no longer restricted by the acoustic stability condition. Written in terms of its Riemann invariants, the pressure subsystem decouples into independent one-dimensional problems along the grid lines, each of them governed by a banded linear system whose bandwidth does not depend on the mesh. The implicit step therefore requires no non-linear solver and costs $O(N_xN_y)$ operations, regardless of the acoustic Courant number. A two-dimensional exactly well-balanced reconstruction, whose limiter weights are frozen within each time step, preserves this linear structure for the second-order scheme. We prove that both the first- and the second-order schemes preserve water-at-rest equilibria exactly, and we assess them on several two-dimensional benchmarks, where they resolve fine-scale velocity structures on coarse meshes and run up to $8.5$ times faster than standard explicit Riemann solvers.}



\maketitle

\section{Introduction}

In this work, we consider the two-dimensional shallow water equations
\begin{equation}\label{eq:shallow-water-2d}
  \left\{ \begin{array}{l}
    \partial_t h +\partial_x \qu + \partial_y \qv  =  0,                           \\
    \partial_t \qu+ \partial_x (\qu u + p) + \partial_y (\qu v)  =  gh\partial_xH, \\
    \partial_t \qv+ \partial_x (\qv u) + \partial_y (\qv v + p)  =  gh\partial_yH,\end{array}\right.
\end{equation}
where $(x,y)\in \Omega \subset \mathbb{R}^2$ represents the spatial domain occupied by the fluid, and $t\in [0,T]$ denotes the time interval. The variable $h$ represents the water height, $\mathbf{u} = (u, v)^T$ is the velocity vector field, $p$ is the hydrostatic pressure, namely $p \equiv p(h) = \tfrac{g}{2}h^2$, with $g$ the gravitational acceleration, and $H$ is a given smooth function representing the bottom topography (the depth, measured from a reference level), which is assumed to be time-independent. Moreover, we denote by $\mathbf{q} = (\qu,\qv)^T= (hu, hv)^T$ the discharge vector field.

System \eqref{eq:shallow-water-2d} can be recast as a balance law in two space dimensions:
\begin{equation}\label{eq:balanced-laws}
  \partial_t \bm{U} + \partial_x \bm{\mathcal{F}}_1(\bm{U}) + \partial_y \bm{\mathcal{F}}_2(\bm{U})  = \bm{\mathcal{S}}_1(\bm{U}) \partial_x H + \bm{\mathcal{S}}_2(\bm{U})\partial_yH,
\end{equation}
where the vector of conserved variables $\bm{U}$, the flux functions $\bm{\mathcal{F}}_1, \bm{\mathcal{F}}_2$, and the source terms $\bm{\mathcal{S}}_1, \bm{\mathcal{S}}_2$, are given by
\begin{equation} \label{eq:2d-bl_components}
  \bm{U} = \begin{pmatrix} h \\ \qu \\ \qv \end{pmatrix}, \
  \bm{\mathcal{F}}_1 = \begin{pmatrix} \qu \\[0mm] \dfrac{\qu^2}{h} + p \\[3mm] \dfrac{\qu \qv}{h}\end{pmatrix}, \
  \bm{\mathcal{F}}_2 = \begin{pmatrix} \qv \\[1mm] \dfrac{\qu \qv}{h} \\[2mm] \dfrac{\qv^2}{h} + p \end{pmatrix}, \
  \bm{\mathcal{S}}_1 = \begin{pmatrix} 0 \\ gh \\ 0 \end{pmatrix}, \
  \bm{\mathcal{S}}_2 = \begin{pmatrix} 0 \\ 0 \\ gh \end{pmatrix}.
\end{equation}

Let $\bm{\mathcal{J}}_i(\bm{U}) = \dfrac{\partial \bm{\mathcal{F}}_i}{\partial \bm{U}}(\bm{U})$ for $i=1,2$ be the Jacobian matrices of the fluxes. For any unit vector $\mathbf{n} = (n_1, n_2)^T \in \mathbb{R}^2$, we define the matrix
\begin{displaymath}
  \bm{\mathcal{A}}(\bm{U},\mathbf{n}) = \bm{\mathcal{J}}_1(\bm{U})n_1 + \bm{\mathcal{J}}_2(\bm{U})n_2.
\end{displaymath}
The two-dimensional system of balance laws \eqref{eq:balanced-laws} is strictly hyperbolic if, for every~$\bm{U}$ in the phase space and every unit vector $\mathbf{n}$, the eigenvalues of $\bm{\mathcal{A}}(\bm{U},\mathbf{n})$ are real and distinct. Under this condition, it admits a factorization
\begin{displaymath}
  \bm{\mathcal{A}}(\bm{U},\mathbf{n}) = \bm{\mathcal{K}}(\bm{U},\mathbf{n}) \; \bm{\mathcal{D}}(\bm{U},\mathbf{n}) \;\bm{\mathcal{K}}^{-1}(\bm{U},\mathbf{n}),
\end{displaymath}
where $\bm{\mathcal{D}}(\bm{U},\mathbf{n})$ is the diagonal matrix containing the eigenvalues of $\bm{\mathcal{A}}(\bm{U},\mathbf{n})$, and $\bm{\mathcal{K}}(\bm{U},\mathbf{n})$ is the matrix whose columns are corresponding right eigenvectors.

System \eqref{eq:balanced-laws}-\eqref{eq:2d-bl_components} is strictly hyperbolic over the physical phase space $\{(h,hu,hv)^T \in \mathbb{R}^3 : h > 0\}$, with eigenvalues given by
\begin{displaymath}
  \lambda_1 = \mathbf{u} \cdot \mathbf{n} - c, \quad \lambda_2 = \mathbf{u} \cdot \mathbf{n},\quad \lambda_3 = \mathbf{u} \cdot \mathbf{n} + c,
\end{displaymath}
where $c = \sqrt{gh}$ represents the celerity of small-amplitude gravity waves.

Depending on the fluid velocity, two distinct flow regimes can be identified: the subcritical regime, where $|\mathbf{u}| < c$, and the supercritical regime, where $|\mathbf{u}| > c$. These regimes are characterized by the Froude number $F_r = \tfrac{|\mathbf{u}|}{c}$, which plays a role analogous to the Mach number in the Euler equations. Specifically, the subcritical regime corresponds to $F_r < 1$, while the supercritical regime corresponds to $F_r > 1$. In this work, we focus primarily on applications within the low-Froude-number regime ($F_r \ll 1$). This is precisely the setting where an implicit treatment of the pressure terms becomes highly advantageous, as it removes the severe time-step restriction imposed by the explicit CFL condition.

An important feature of the shallow water equations is the existence of steady-state solutions, or equilibria, satisfying
\begin{equation*}
  \left\{ \begin{array}{l}
    \partial_x \qu + \partial_y \qv = 0,                        \\
    \partial_x (\qu u+p) + \partial_y (\qu v) = gh\partial_x H, \\
    \partial_x (\qv u) + \partial_y (\qv v+p) = gh\partial_y H.\end{array}\right.
\end{equation*}
In this work, we are interested in developing a scheme that exactly preserves the family of water-at-rest stationary solutions, defined by
\begin{equation}\label{eq:water-at-rest-shallow-water}
  \left\{ \begin{array}{l}
    u = v = 0,                     \\
    \partial_x p = gh\partial_x H, \\
    \partial_y p = gh\partial_y H.\end{array}\right.
\end{equation}
Numerical schemes capable of preserving this specific class of equilibria are referred to as exactly well-balanced schemes, since they exactly preserve the steady states and not just a discrete approximation of them.

Consequently, our goal is to design a semi-implicit, exactly well-balanced scheme, focusing on computational efficiency, performance in low-Froude-number regimes and robust applicability to real-world scenarios. The splitting strategy employed here can be motivated by a suitable rescaling of the governing equations. Indeed, in the dimensionless form of the shallow water equations, the scaling factors amplify the role of the pressure terms as the Froude number approaches zero. To address the resulting multi-scale nature of the flow, we combine the Suliciu relaxation technique introduced in~\cite{bouchut2004nonlinear} with the convective flux-splitting strategy inspired by \cite{franck2020semi} for the Ripa equations.

The design of numerical methodologies capable of accurately maintaining balanced states has been the subject of extensive study in many different works, including \cite{bermudez1994upwind, greenberg1996well, perthame2001kinetic, gallice2002solveurs, audusse2004fast, chinnayya2004well, berthon2012efficient, audusse2015simple}. In this context, our primary interest lies in preserving the unperturbed water-at-rest steady states given by \eqref{eq:water-at-rest-shallow-water}, meaning that the scheme must exactly preserve states where the velocity field satisfies $\mathbf{u} = \mathbf{0}$.

At the same time, different techniques for the design of implicit or semi-implicit schemes for shallow water models have been continuously developed since the pioneering work \cite{casulli1990semi}. For instance, implicit-explicit (IMEX) methods were proposed for three‐dimensional shallow water flows in \cite{casulli1992semi, casulli1994stability, casulli1997numerical, casulli2000unstructured}, while more general IMEX frameworks have been successfully applied to hyperbolic systems in works such as \cite{boscarino2009class, boscarino2017unified}. These implicit and semi-implicit approaches have been explored using a wide variety of spatial discretizations, including finite volume methods \cite{bispen2014imex, vater2018semi, cea2020extension}, discontinuous Galerkin (DG) formulations \cite{giraldo2009high, dumbser2013staggered, tumolo2013semi, tavelli2014high, kang2020imex}, finite difference approximations \cite{casulli1990semi, casulli1992semi}, and hybrid or staggered strategies \cite{bonaventura2018multilayer, busto2022staggered, boscheri2023all}.

In general, the core idea behind these multi-scale methods is to separate the fast acoustic waves from the slow transport ones, and then to couple explicit and implicit time advances. Particularly in the low-Froude-number limit, a notable IMEX finite volume approach was proposed in \cite{bispen2014imex}. More recently, exactly well-balanced implicit and IMEX Lagrange-projection finite volume schemes have been developed for the 1D shallow water system in \cite{caballero2023implicit}. This framework was later extended to a semi-implicit fully exactly well-balanced relaxation scheme for the 1D shallow water equations in~\cite{caballero2024semi} and to the two-layer shallow water system in \cite{caballero2025semi}.

Unlike traditional IMEX formulations, which frequently suffer from the computational bottleneck of solving complex non-linear systems, our framework substantially simplifies the implicit step. The proposed splitting strategy yields a pressure subsystem that reduces to a linear system at each time step.

Although the system matrix is dynamically updated to reflect changing states, its specific structure ensures an efficient resolution of the system. This approach entirely bypasses the need for costly iterative non-linear solvers, thereby minimizing computational overhead and yielding a highly efficient, scalable framework for large-scale simulations.

The resulting linear system is solved directly, with a cost that is independent of the CFL number associated with the pressure step.

To the best of our knowledge, the integration of an exactly well-balanced property within a semi-implicit relaxation splitting framework has not yet been extended to the two-dimensional shallow water equations \eqref{eq:shallow-water-2d}. The present work introduces the following contributions.

First, we provide a rigorous treatment of the tangential terms arising in the two-dimensional formulation, demonstrating that these contributions can be neglected while preserving the well-balanced property exactly up to machine precision. This result yields a directional decoupling of the pressure system, substantially reducing its computational complexity without compromising the scheme's structural properties.

Second, we develop a exactly well-balanced two-dimensional reconstruction procedure equipped with a limiter that preserves the linearity of the pressure system upon linearization, thereby enabling the extension of the method to second-order accuracy while retaining the computational advantages of the underlying splitting.

Third, we show that the resulting scheme performs well in capturing multidimensional velocity structures, remaining robust in the presence of sharp discontinuities.

The remainder of this paper is organized as follows. Section~\ref{sec:split-and-relax} introduces the relaxation framework for the 2D shallow water equations and details the splitting of the relaxed system into pressure and transport subsystems. Section~\ref{sec:semi-discrete} presents the semi-discrete formulation of the scheme, which is subsequently enhanced to exactly satisfy the well-balanced property in Section~\ref{sec:well-balanced}. The corresponding first- and second-order  exactly well-balanced reconstruction operators are described in Section~\ref{sec:well-balance-reconstruction-operators}. Section~\ref{sec:fully-discrete} defines the fully discrete version of the proposed methodology. Finally, Section~\ref{sec:numerical-experiments} reports several numerical experiments to assess the accuracy, robustness, and computational efficiency of the proposed schemes.

\section{Splitting of the relaxation scheme}\label{sec:split-and-relax}

We begin by proposing a relaxed system for (\ref{eq:shallow-water-2d}), employing relaxation techniques inspired by the framework introduced in~\cite{suliciu1990modelling,suliciu1992some,bouchut2004nonlinear}. This system is formulated as follows:
\begin{equation}\label{eq:full-relaxed-system}
  \left\{ \begin{array}{l}
    \partial_t h +\partial_x \qu + \partial_y \qv  =  0,                               \\
    \partial_t \qu + \partial_x (\qu u + \pi) + \partial_y (\qu v)  =  gh\partial_x H, \\
    \partial_t \qv + \partial_x (\qv u) + \partial_y (\qv v + \pi)  =  gh\partial_y H, \\
    \partial_t \Pi + \partial_x(\Pi u + a^2 u)  + \partial_y(\Pi v + a^2 v) = -\frac{1}{\varepsilon} (\pi - p),
  \end{array}\right.
\end{equation}
with $\Pi = h \pi$ and where $a$ is a constant satisfying the subcharacteristic condition, that is,
\begin{displaymath}
  a \geq h \sqrt{gh}.
\end{displaymath}

Note that in the limit $\varepsilon \rightarrow 0$, we formally recover $\pi = p$. Consequently, the $\varepsilon$-source term in (\ref{eq:full-relaxed-system}) is typically omitted and the relaxation is applied instantaneously by setting $\pi$ equal to the hydrostatic pressure $p=\tfrac{g}{2}h^2$ at each time step.

The proposed relaxed system can be expressed as
\begin{equation}\label{eq:relaxed-balanced-laws}
  \partial_t \bm{W} + \partial_x {\bm{F}}_1(\bm{W}) + \partial_y {\bm{F}}_2(\bm{W}) + \partial_x {\bm{G}}_1(\bm{W}) + \partial_y {\bm{G}}_2(\bm{W}) = {\bm{S}}_1(\bm{W}) \partial_x H + {\bm{S}}_2(\bm{W})\partial_yH,
\end{equation}
where
\begin{equation}\label{eq:vector-w}
  \bm{W} = \begin{pmatrix} h \\ \qu \\ \qv \\ \Pi \end{pmatrix},
\end{equation}
the fluxes are given by
\begin{equation}\label{eq:relaxed-fluxes}
  \bm{F}_{1} = \begin{pmatrix} \qu \\[1mm] \frac{\qu^2}{h} \\[1mm] \frac{\qu \qv}{h} \\[1mm] \frac{\Pi \qu}{h} \end{pmatrix},\
  \bm{F}_{2} = \begin{pmatrix} \qv \\[1mm] \frac{\qu \qv}{h} \\[1mm] \frac{\qv^2}{h} \\[1mm] \frac{\Pi \qv}{h} \end{pmatrix},\
  \bm{G}_{1} = \begin{pmatrix} 0 \\[1mm] \frac{\Pi}{h} \\[1mm] 0 \\[1mm] a^2 \dfrac{\qu}{h}\end{pmatrix},\
  \bm{G}_{2} = \begin{pmatrix} 0 \\[1mm] 0 \\[1mm] \frac{\Pi}{h} \\[1mm] a^2 \dfrac{\qv}{h}\end{pmatrix},
\end{equation}
and
\begin{equation}\label{eq:relaxed-sources}
  {\bm{S}}_1 = \begin{pmatrix} 0 \\ gh \\ 0 \\ 0 \end{pmatrix},\
  {\bm{S}}_2 = \begin{pmatrix} 0 \\ 0 \\ gh \\ 0 \end{pmatrix}.
\end{equation}

This formulation can then be decomposed into two distinct subsystems by a splitting technique:
\begin{equation}\label{eq:balance-law-pressure}
  \partial_t \bm{W} +  \partial_x {\bm{G}}_1(\bm{W}) + \partial_y {\bm{G}}_2(\bm{W}) = {\bm{S}}_1(\bm{W}) \partial_x H + {\bm{S}}_2(\bm{W})\partial_yH,
\end{equation}
and
\begin{equation}\label{eq:balance-law-transport}
  \partial_t \bm{W} + \partial_x {\bm{F}}_1(\bm{W}) + \partial_y {\bm{F}}_2(\bm{W}) = 0,
\end{equation}
to be solved sequentially at every time step. System (\ref{eq:balance-law-pressure}) is referred to as the pressure system, since it contains the pressure terms, and reads
\begin{equation}\label{eq:pressure-system}
  \left\{ \begin{array}{l}
    \partial_t h = 0,                                   \\
    \partial_t \qu + \partial_x  \pi =  gh\partial_x H, \\
    \partial_t \qv + \partial_y  \pi =  gh\partial_y H, \\[0.5mm]
    \partial_t \Pi + a^2 \partial_x u + a^2 \partial_y v = 0.
  \end{array}\right.
\end{equation}
Conversely, system (\ref{eq:balance-law-transport}) is identified as the transport system, and reads
\begin{equation}\label{eq:transport-system}
  \left\{ \begin{array}{l}
    \partial_t h + \partial_x (hu) + \partial_y (hv) = 0,        \\
    \partial_t \qu + \partial_x (\qu u) + \partial_y(\qu v) = 0, \\
    \partial_t \qv + \partial_x (\qv u) + \partial_y(\qv v) = 0, \\
    \partial_t \Pi + \partial_x (\Pi u) + \partial_y (\Pi v) = 0
    .\end{array}\right.
\end{equation}

It can be verified that the eigenvalues of system \eqref{eq:pressure-system} are
\begin{equation}\label{eq:eigenvalues_pressure}
  \lambda_- = -\dfrac{a}{h}, \quad \lambda_0 = 0 \text{ (double)}, \quad \lambda_+ = \dfrac{a}{h},
\end{equation}
while system \eqref{eq:transport-system} has eigenvalues given by
\begin{equation}\label{eq:eigenvalues_transport}
  \lambda = \mathbf{u} \cdot \mathbf{n} \text{ (quadruple)}.
\end{equation}

The projection of the pressure system \eqref{eq:pressure-system} along any normal direction can be rewritten in terms of its Riemann invariants, defined as
\begin{equation}\label{eq:riemann-invariants}
  \wR = \pi + a \ \mathbf{u} \cdot \mathbf{n}, \quad \wL = \pi - a \ \mathbf{u} \cdot \mathbf{n},
\end{equation}
giving the following differential system
\begin{equation}\label{eq:riemann-invariants-system}
  \left\{ \begin{array}{rcl}
    \partial_t \wL & = & \dfrac{a}{h} \nabla \wL \cdot \mathbf{n}- \dfrac{a^2}{h}  \nabla(\mathbf{u} \cdot \mathbf{t}) \cdot \mathbf{t} - ag \nabla H \cdot \mathbf{n},   \\
    \partial_t \wR & = & -\dfrac{a}{h} \nabla \wR \cdot \mathbf{n} - \dfrac{a^2}{h} \nabla (\mathbf{u} \cdot \mathbf{t}) \cdot \mathbf{t} + ag \nabla H \cdot \mathbf{n},
  \end{array}\right.
\end{equation}
where $\mathbf{t} = (-n_2,n_1)^T$ represents the orthogonal unit tangent vector, and the water height $h$ is treated as frozen in time, consistent with the first equation of \eqref{eq:pressure-system}—in the sense detailed later in Section \ref{sec:well-balanced}.

The primary advantage of this splitting approach lies in the decoupling of acoustic and transport phenomena. The pressure system can be integrated either explicitly or implicitly, whereas the transport system is always solved using an explicit scheme. At low Froude numbers, the time-step constraint is predominantly governed by the pressure term. Hence, treating the pressure step implicitly allows for longer time steps, since it leads to a significantly less restrictive CFL condition.

\section{Semi-discrete schemes of the subsystems}\label{sec:semi-discrete}

\begin{figure}[htbp]
  \centering
  \includegraphics{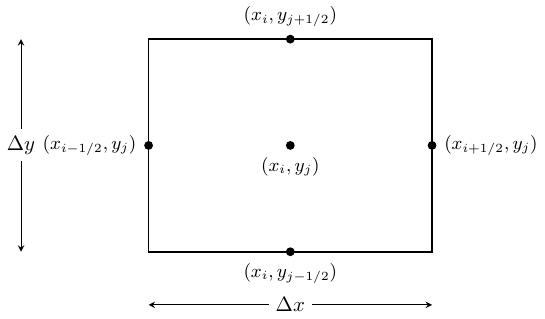}
  \caption{Schematic representation of the Cartesian control volume $C_{i,j}$ used for the 2D finite volume discretization.}
  \label{fig:fv_cell}
\end{figure}

We consider a structured discretization of the computational domain $\Omega = [a,b]\times[c,d]$ into uniform control cells $C_{i,j}$ of constant volume $|C_{i,j}| = \Delta x \, \Delta y$, where $\Delta x$ and $\Delta y$ denote the constant spatial steps. The cell interfaces are defined by $x_{i-1/2} = a + i \Delta x$ and $y_{j-1/2} = c + j \Delta y$ for $i \in \{0,\dots,N_x-1\}$ and $j \in \{0,\dots,N_y-1\}$. Consequently, the control cells are given by
\begin{equation*}
  C_{i,j} = [x_{i-1/2}, x_{i+1/2}] \times [y_{j-1/2}, y_{j+1/2}]
\end{equation*}
while cell centres are given by
\begin{equation*}
  (x_i, y_j) = \left(\tfrac{x_{i-1/2}+x_{i+1/2}}{2}, \tfrac{y_{j-1/2}+y_{j+1/2}}{2}\right),
\end{equation*}
see Figure \ref{fig:fv_cell}.
The time variable will remain continuous for now, and the cell average approximations are denoted by
\begin{displaymath}
  \bm{W}_{i,j}(t) \approx \dfrac{1}{|C_{i,j}|} \int_{C_{i,j}} \bm{W}(x,y,t)  \; dx \; dy,
\end{displaymath}
with $\bm{W}$ denoting the corresponding vector of unknowns, extending conserved variables as
\begin{equation*}
  \bm{W} = \begin{pmatrix} \bm{U} \\ \Pi \end{pmatrix}.
\end{equation*}
A semi-discrete finite volume scheme for pressure system~(\ref{eq:pressure-system}) and transport system~(\ref{eq:transport-system}) can be written as
\begin{equation}\label{eq:semi-discrete-pressure-system}
  \left\{ \begin{array}{rcl}
    h^\prime_{i,j}(t)       & = & 0,                                                                                                \\[1mm]

    {(\qu)}^\prime_{i,j}(t) & = & - \dfrac{1}{\Delta x} \Big({\pi}^{*}_{i+1/2,j}(t) - {\pi}^{*}_{i-1/2,j}(t)\Big)+({S}_1)_{i,j}(t), \\[3mm]

    {(\qv)}^\prime_{i,j}(t) & = & - \dfrac{1}{\Delta y} \Big({\pi}^{*}_{i,j+1/2}(t) - {\pi}^{*}_{i,j-1/2}(t)\Big)+({S}_2)_{i,j}(t), \\[3mm]

    \Pi^\prime_{i,j}(t)     & = & - \dfrac{a^2}{\Delta x} \Big({u}^{*}_{i+1/2,j}(t) - {u}^{*}_{i-1/2,j}(t)\Big)                     \\[3mm]
                            &   & - \dfrac{a^2}{\Delta y} \Big({v}^{*}_{i,j+1/2}(t)- {v}^{*}_{i,j-1/2}(t)\Big),
  \end{array}\right.
\end{equation}
where
\begin{displaymath}
  \begin{array}{rcl}
    ({S}_1)_{i,j}(t) & \approx & \dfrac{1}{|C_{i,j}|} \displaystyle\int_{C_{i,j}} gh(x,y,t) \partial_x H(x,y) \;dx\;dy,   \\[3mm]
    ({S}_2)_{i,j}(t) & \approx & \dfrac{1}{|C_{i,j}|} \displaystyle\int_{C_{i,j}} gh(x,y,t) \,\partial_y H(x,y) \;dx\;dy, \\
  \end{array}
\end{displaymath}
and
\begin{equation}\label{eq:semi-discrete-transport-system}
  \left\{ \begin{array}{rcl}
    h^\prime_{i,j}(t)       & = & - \dfrac{1}{\Delta x} \Big(h^{*}_{i+1/2,j}(t){u}^{*}_{i+1/2,j}(t)-h^{*}_{i-1/2,j}(t){u}^{*}_{i-1/2,j}(t) \Big)             \\[3mm]
                            &   & - \dfrac{1}{\Delta y} \Big(h^{*}_{i,j+1/2}(t)v^{*}_{i,j+1/2}(t) - h^{*}_{i,j-1/2}(t)v^{*}_{i,j-1/2}(t)\Big),               \\[3mm]

    {(\qu)}^\prime_{i,j}(t) & = & - \dfrac{1}{\Delta x} \Big({(\qu)}^{*}_{i+1/2,j}(t){u}^{*}_{i+1/2,j}(t)-{(\qu)}^{*}_{i-1/2,j}(t){u}^{*}_{i-1/2,j}(t) \Big) \\[3mm]
                            &   & -\dfrac{1}{\Delta y} \Big({(\qu)}^{*}_{i,j+1/2}(t){v}^{*}_{i,j+1/2}(t)-{(\qu)}^{*}_{i,j-1/2}(t){v}^{*}_{i,j-1/2}(t) \Big), \\[3mm]

    {(\qv)}^\prime_{i,j}(t) & = & - \dfrac{1}{\Delta x} \Big({(\qv)}^{*}_{i+1/2,j}(t){u}^{*}_{i+1/2,j}(t)-{(\qv)}^{*}_{i-1/2,j}(t){u}^{*}_{i-1/2,j}(t) \Big) \\[3mm]
                            &   & -\dfrac{1}{\Delta y} \Big({(\qv)}^{*}_{i,j+1/2}(t){v}^{*}_{i,j+1/2}(t)-{(\qv)}^{*}_{i,j-1/2}(t){v}^{*}_{i,j-1/2}(t) \Big), \\[3mm]

    \Pi^\prime_{i,j}(t)     & = & - \dfrac{1}{\Delta x} \Big({(\Pi)}^{*}_{i+1/2,j}(t){u}^{*}_{i+1/2,j}(t)-{(\Pi)}^{*}_{i-1/2,j}(t){u}^{*}_{i-1/2,j}(t) \Big) \\[3mm]
                            &   & -\dfrac{1}{\Delta y} \Big({(\Pi)}^{*}_{i,j+1/2}(t){v}^{*}_{i,j+1/2}(t)-{(\Pi)}^{*}_{i,j-1/2}(t){v}^{*}_{i,j-1/2}(t) \Big),
  \end{array}\right.
\end{equation}
respectively. In the previous expressions, $X_{i\pm1/2,j}^*$, $X_{i,j\pm1/2}^*$ denote a suitable approximation of $X$ at the centre of cell edges, that is
\begin{equation*}
  X_{i\pm1/2,j}^*(t)\approx X(x_{i\pm 1/2},y_j,t), \quad X_{i,j\pm1/2}^*(t)\approx X(x_i,y_{j\pm 1/2},t),
\end{equation*}
for $X\in \{\pi,u,v\}$, and also $X\in \{h,q_u,q_v,\Pi\}$. These approximations will be detailed further below. However, a scheme like the previous one will not be well-balanced in general.

\section{Well-balanced schemes}\label{sec:well-balanced}

To derive well-balanced schemes preserving the water-at-rest states \eqref{eq:water-at-rest-shallow-water} exactly, we adopt the general methodology outlined in \cite{castro2020well}. These steady states are characterized by zero velocity. Consequently, the well-balanced correction does not affect the transport system, which already preserves states at rest, and needs to be applied only to the pressure system. Moreover, the approach followed here is strictly local. Over a given time step $[t^n, t^{n+1}]$ (with $t^{n+1} = t^n + \Delta t$), we consider a stationary solution within each control cell $C_{i,j}$ that satisfies \eqref{eq:water-at-rest-shallow-water} at time $t=t^n$. Specifically, we define $\bm{U}_{i,j}^{s,n} = (h^{s,n}_{i,j}, 0, 0)^T$ where
\begin{equation*}
  h^{s,n}_{i,j}(x,y) =  h_{i,j}(t^n) - H(x_i,y_j) + H(x,y).
\end{equation*}
Following the finite volume approach, we approximate the second and third equations of \eqref{eq:water-at-rest-shallow-water} within the control cell $C_{i,j}$ using the mid-point rule to compute the integral, obtaining
\begin{equation}\label{eq:steady-states-semi-discrete}
  \begin{aligned}
    \frac{1}{\Delta x} \Big(\pi^{s,n}_{i,j}(x_{i+1/2}, y_j) - \pi^{s,n}_{i,j}(x_{i-1/2}, y_j)\Big) 
    \approx {} \frac{1}{|C_{i,j}|} \int_{C_{i,j}} gh_{i,j}^{s,n}(x,y) \,\partial_x H(x,y) \,dxdy,                            & \\[1ex]
    \frac{1}{\Delta y} \Big(\pi^{s,n}_{i,j}(x_{i}, y_{j+1/2}) - \pi^{s,n}_{i,j}(x_{i}, y_{j-1/2})\Big)  
    \approx {} \frac{1}{|C_{i,j}|} \int_{C_{i,j}} gh_{i,j}^{s,n}(x,y) \,\partial_y H(x,y) \,dxdy,                            &
  \end{aligned}
\end{equation}
where $\pi^{s,n}_{i,j}$ reduces to the mean hydrostatic pressure in the cell, that is $\pi^{s,n}_{i,j}=p^{s,n}_{i,j}=\tfrac{1}{2}g(h^{s,n}_{i,j})^2$. Adding the relations in \eqref{eq:steady-states-semi-discrete} to the second and third equations of the pressure system \eqref{eq:semi-discrete-pressure-system} yields:
\begin{equation}\label{eq:semi-discrete-pressure-system-well-balanced}
  \left\{ \begin{array}{rcl}
    h^\prime_{i,j}(t)       & = & 0,                                                                                                                      \\[1mm]

    {(\qu)}^\prime_{i,j}(t) & = & - \dfrac{1}{\Delta x} \Big({\pi}^{*}_{i+1/2,j}(t) - \pi^{s,n}_{i,j}(x_{i+1/2}, y_j)\Big)                                \\[3mm]
                            &   & + \dfrac{1}{\Delta x} \Big({\pi}^{*}_{i-1/2,j}(t) - \pi^{s,n}_{i,j}(x_{i-1/2}, y_j)\Big) + ({\mathcal{S}}_1)_{i,j}(t),  \\[3mm]

    {(\qv)}^\prime_{i,j}(t) & = & - \dfrac{1}{\Delta y} \Big({\pi}^{*}_{i,j+1/2}(t) - \pi^{s,n}_{i,j}(x_{i}, y_{j+1/2})\Big)                              \\[3mm]
                            &   & + \dfrac{1}{\Delta y} \Big({\pi}^{*}_{i,j-1/2}(t) - \pi^{s,n}_{i,j}(x_{i}, y_{j-1/2})\Big)+ ({\mathcal{S}}_2)_{i,j}(t), \\[3mm]

    \Pi^\prime_{i,j}(t)     & = & - \dfrac{a^2}{\Delta x} \Big({u}^{*}_{i+1/2,j}(t) - {u}^{*}_{i-1/2,j}(t)\Big)                                           \\[3mm]
                            &   & - \dfrac{a^2}{\Delta y} \Big({v}^{*}_{i,j+1/2}(t)- {v}^{*}_{i,j-1/2}(t)\Big),
  \end{array}\right.
\end{equation}
where
\begin{displaymath}
  \begin{array}{rcl}
    ({\mathcal{S}}_1)_{i,j}(t) & \approx & \dfrac{1}{|C_{i,j}|} \displaystyle\int_{C_{i,j}} g(h(x,y,t) - h^{s,n}_{i,j}(x,y))\partial_x H(x,y) \;dx\;dy, \\[3mm]
    ({\mathcal{S}}_2)_{i,j}(t) & \approx & \dfrac{1}{|C_{i,j}|} \displaystyle\int_{C_{i,j}} g(h(x,y,t) - h^{s,n}_{i,j}(x,y))\partial_y H(x,y) \;dx\;dy. \\
  \end{array}
\end{displaymath}

Furthermore, the first equation in \eqref{eq:semi-discrete-pressure-system-well-balanced} implies that $h_{i,j}$ remains constant throughout the time step $[t^n, t^{n+1}]$, frozen at its initial value at $t=t^n$. We can then rewrite the pressure system \eqref{eq:semi-discrete-pressure-system-well-balanced} as follows:
\begin{equation}\label{eq:semi-discrete-pressure-system-wb}
  \left\{ \begin{array}{rcl}

    {u}^\prime_{i,j}(t) & = & - \dfrac{1}{h_{i,j}(t^n) \Delta x} \Big({\pi}^{*}_{i+1/2,j}(t) - \pi^{s,n}_{i,j}(x_{i+1/2}, y_j)\Big)   \\[3mm]
                        &   & + \dfrac{1}{h_{i,j}(t^n)\Delta x} \Big({\pi}^{*}_{i-1/2,j}(t) - \pi^{s,n}_{i,j}(x_{i-1/2}, y_j)\Big)    \\[3mm]
                        &   & + \dfrac{1}{h_{i,j}(t^n) } ({\mathcal{S}}_1)_{i,j}^{n},                                                 \\[3mm]

    {v}^\prime_{i,j}(t) & = & - \dfrac{1}{h_{i,j}(t^n) \Delta y} \Big({\pi}^{*}_{i,j+1/2}(t) - \pi^{s,n}_{i,j}(x_{i}, y_{j+1/2})\Big) \\[3mm]
                        &   & + \dfrac{1}{h_{i,j}(t^n) \Delta y} \Big({\pi}^{*}_{i,j-1/2}(t) - \pi^{s,n}_{i,j}(x_{i}, y_{j-1/2})\Big) \\[3mm]
                        &   & + \dfrac{1}{h_{i,j}(t^n)} ({\mathcal{S}}_2)_{i,j}^{n},                                                  \\[3mm]

    \pi^\prime_{i,j}(t) & = & - \dfrac{a^2}{h_{i,j}(t^n) \Delta x} \Big({u}^{*}_{i+1/2,j}(t) - {u}^{*}_{i-1/2,j}(t) \Big)             \\[3mm]
                        &   & - \dfrac{a^2}{h_{i,j}(t^n) \Delta y} \Big({v}^{*}_{i,j+1/2}(t)- {v}^{*}_{i,j-1/2}(t) \Big),
  \end{array}\right.
\end{equation}
with
\begin{equation*}
  \begin{array}{rcl}
    ({\mathcal{S}}_1)_{i,j}^{n} & \approx & \dfrac{1}{|C_{i,j}|} \displaystyle\int_{C_{i,j}} g(h(x,y,t^n) - h^{s,n}_{i,j}(x,y))\partial_x H(x,y) \;dx\;dy, \\[3mm]
    ({\mathcal{S}}_2)_{i,j}^{n} & \approx & \dfrac{1}{|C_{i,j}|} \displaystyle\int_{C_{i,j}} g(h(x,y,t^n) - h^{s,n}_{i,j}(x,y))\partial_y H(x,y) \;dx\;dy.
  \end{array}
\end{equation*}
System \eqref{eq:semi-discrete-pressure-system-wb} is rewritten in terms of the Riemann invariants of the pressure system \eqref{eq:riemann-invariants} along the principal directions as
\begin{equation*}
  \wRx = \pi + a \ u, \quad \wLx = \pi - a \ u,
\end{equation*}
in the horizontal direction and
\begin{equation*}
  \wRy = \pi + a \ v, \quad \wLy = \pi - a \ v,
\end{equation*}
in the vertical direction. This yields the following semi-discrete system in the horizontal direction
\begin{equation}\label{eq:semi-discrete-pressure-system-wb-horizontal}
  \left\{ \begin{array}{rcl}
    {(\wLx)}^\prime_{i,j}(t) & = & + \dfrac{a}{h_{i,j}(t^n) \Delta x} \Big({(\wLx)}_{i+1/2,j}(t) - (\wLx)^{s,n}_{i,j}(x_{i+1/2}, y_j)\Big) \\[3mm]
                             &   & - \dfrac{a}{h_{i,j}(t^n) \Delta x} \Big({(\wLx)}_{i-1/2,j}(t) - (\wLx)^{s,n}_{i,j}(x_{i-1/2}, y_j)\Big) \\[3mm]
                             &   & -\dfrac{a^2}{h_{i,j}(t^n) \Delta y} \Big({v}^{*}_{i,j+1/2}(t)- {v}^{*}_{i,j-1/2}(t) \Big)               \\[3mm]
                             &   & - \dfrac{a}{h_{i,j}(t^n)} ({\mathcal{S}}_1)_{i,j}^{n},                                                  \\[3mm]

    {(\wRx)}^\prime_{i,j}(t) & = & - \dfrac{a}{h_{i,j}(t^n) \Delta x} \Big({(\wRx)}_{i+1/2,j}(t) - (\wRx)^{s,n}_{i,j}(x_{i+1/2}, y_j)\Big) \\[3mm]
                             &   & + \dfrac{a}{h_{i,j}(t^n) \Delta x} \Big({(\wRx)}_{i-1/2,j}(t) - (\wRx)^{s,n}_{i,j}(x_{i-1/2}, y_j)\Big) \\[3mm]
                             &   & - \dfrac{a^2}{h_{i,j}(t^n) \Delta y} \Big({v}^{*}_{i,j+1/2}(t)- {v}^{*}_{i,j-1/2}(t) \Big)              \\[3mm]
                             &   & + \dfrac{a}{h_{i,j}(t^n)} ({\mathcal{S}}_1)_{i,j}^{n},                                                  \\[3mm]
  \end{array}\right.
\end{equation}
and in the vertical one
\begin{equation}\label{eq:semi-discrete-pressure-system-wb-vertical}
  \left\{ \begin{array}{rcl}
    {(\wLy)}^\prime_{i,j}(t) & = & + \dfrac{a}{h_{i,j}(t^n) \Delta y} \Big({(\wLy)}_{i,j+1/2}(t) - (\wLy)^{s,n}_{i,j}(x_{i}, y_{j+1/2})\Big) \\[3mm]
                             &   & - \dfrac{a}{h_{i,j}(t^n) \Delta y} \Big({(\wLy)}_{i,j-1/2}(t) - (\wLy)^{s,n}_{i,j}(x_{i}, y_{j-1/2})\Big) \\[3mm]
                             &   & - \dfrac{a^2}{h_{i,j}(t^n) \Delta x} \Big({u}^{*}_{i+1/2,j}(t)- {u}^{*}_{i-1/2,j}(t)) \Big)               \\[3mm]
                             &   & - \dfrac{a}{h_{i,j}(t^n)} ({\mathcal{S}}_2)_{i,j}^{n},                                                    \\[3mm]

    {(\wRy)}^\prime_{i,j}(t) & = & - \dfrac{a}{h_{i,j}(t^n) \Delta y} \Big({(\wRy)}_{i,j+1/2}(t) - (\wRy)^{s,n}_{i,j}(x_{i}, y_{j+1/2})\Big) \\[3mm]
                             &   & + \dfrac{a}{h_{i,j}(t^n) \Delta y} \Big({(\wRy)}_{i,j-1/2}(t) - (\wRy)^{s,n}_{i,j}(x_{i}, y_{j-1/2})\Big) \\[3mm]
                             &   & - \dfrac{a^2}{h_{i,j}(t^n) \Delta x} \Big({u}^{*}_{i+1/2,j}(t)- {u}^{*}_{i-1/2,j}(t)) \Big)               \\[3mm]
                             &   & + \dfrac{a}{h_{i,j}(t^n)} ({\mathcal{S}}_2)_{i,j}^{n}.                                                    \\[3mm]
  \end{array}\right.
\end{equation}
\begin{remark}
  Note that systems \eqref{eq:semi-discrete-pressure-system-wb-horizontal} and \eqref{eq:semi-discrete-pressure-system-wb-vertical} are discretizations of \eqref{eq:riemann-invariants-system}. Moreover, in practice we can use the Riemann invariants to compute the values $X^*$ for $X\in \{\pi, u, v\}$ in \eqref{eq:semi-discrete-pressure-system} and \eqref{eq:semi-discrete-transport-system} as detailed below.
\end{remark}

\section{Well-balanced reconstruction operators}\label{sec:well-balance-reconstruction-operators}
The core components of the approach introduced in \cite{castro2020well} include a exactly well-balanced reconstruction operator, an appropriate quadrature formula, and a suitable approximation of the source term.

A exactly well-balanced reconstruction operator is defined as follows:
\begin{definition} Given a stationary solution $\bm{U}^s$ of a system of balance laws~(\ref{eq:balanced-laws}), a reconstruction operator is said to be exactly well-balanced for $\bm{U}^s$ if, for every pair of indices $\{i,j\}$,
  \begin{displaymath}
    \bm{P}_{i,j}(x,y) = \bm{U}^s(x,y), \quad \forall (x,y) \in C_{i,j},
  \end{displaymath}
  where $\bm{P}_{i,j}(x,y) =\bm{P}_{i,j}(x,y;\{\bm{U}^s_{k,l}\}_{\{k,l\}\in {\mathcal{T}}_{C_{i,j}}})$ and ${\mathcal{T}}_{C_{i,j}}$ denotes the stencil of the reconstruction operator.
\end{definition}

Since we are interested in first- and second-order numerical schemes, the mid-point rule is employed as the quadrature formula. Consequently, as previously stated, we identify cell averages with point values at the cell centres, an approximation that holds up to second-order accuracy. We determine the steady state satisfying (\ref{eq:water-at-rest-shallow-water}) such that its centred point values coincide with the cell averages at the beginning of each time step ($t = t^n$), establishing the condition
\begin{equation}\label{eq:wb-compatibility-condition}
  \bm{U}^{s,n}_{i,j}(x_i, y_j) = \bm{U}^n_{i,j}
\end{equation}
with $\bm{U}^n_{i,j} = \bm{U}_{i,j}(t^n)$.

For the sake of completeness, we now detail the full procedure introduced in \cite{castro2020well} to obtain exactly well-balanced reconstruction operators over the time interval $[t^n, t^{n+1}]$. Note that this process is also designed to preserve a family of steady-state solutions for a system of balance laws.
\begin{enumerate}
  \item At time level $t^n$, identify a stationary solution $\bm{U}^{s,n}_{i,j}$, if one exists, within the stencil of cell $C_{i,j}$ that satisfies
        \begin{displaymath}
          \dfrac{1}{|C_{i,j}|} \int_{C_{i,j}} \bm{U}^{s,n}_{i,j}(x,y) \; dx\;dy = \bm{U}^n_{i,j}.
        \end{displaymath}
        Otherwise, set $\bm{U}^{s,n}_{i,j} \equiv 0$.

  \item Compute the spatial fluctuations at $t^n$ within the cell stencil given by
        \begin{equation}\label{eq:spatial-fluctuation}
          \bm{V}^n_{k,l} = \bm{U}^n_{k,l} - \dfrac{1}{|C_{k,l}|} \int_{C_{k,l}} \bm{U}^{s,n}_{i,j}(x,y) \; dx\;dy, \quad \{k,l\} \in {\mathcal{T}}_{C_{i,j}}.
        \end{equation}
        Subsequently, apply a standard reconstruction operator of the desired order~$p$ (in our case $p=1$ or $p=2$), denoted by $\bm{Q}$, to these fluctuations:
        \begin{displaymath}
          \bm{Q}^n_{i,j}(x,y) = \bm{Q}_{i,j}(x,y; \{\bm{V}^n_{k,l}\}_{\{k,l\} \in {\mathcal{T}}_{C_{i,j}}}).
        \end{displaymath}

  \item Finally, define the spatial exactly well-balanced reconstruction operator as
        \begin{displaymath}
          \bm{P}^n_{i,j}(x,y) = \bm{U}^{s,n}_{i,j}(x,y) + \bm{Q}^n_{i,j}(x,y).
        \end{displaymath}

\end{enumerate}

Nevertheless, these reconstructions do not inherently account for temporal evolution. To resolve this, the well-balanced reconstruction operator at time $t \in [t^n, t^{n+1}]$ is defined as the sum of the well-balanced spatial reconstruction at $t^n$ and a standard reconstruction of the temporal fluctuations. The resulting operator is formulated as
\begin{displaymath}
  \bm{P}_{i,j}(x,y, t) = \bm{P}^n_{i,j}(x,y) + \tilde{\bm{Q}}_{i,j}(x,y,t), \quad (x,y) \in C_{i,j}, \   t \in [t^n, t^n + \Delta t],
\end{displaymath}
where $\tilde{\bm{Q}}_{i,j}(x,y,t)$ is a reconstruction operator defined in terms of time fluctuations
\begin{equation}\label{eq:time-fluctuation}
  \tilde{\bm{V}}^n_{k,l}(t) = \bm{U}_{k,l}(t) - \bm{U}^n_{k,l} \quad \forall \{k,l\} \in {\mathcal{T}}_{C_{i,j}},
\end{equation}
that is,
\begin{displaymath}
  \tilde{\bm{Q}}_{i,j}(x,y,t) = \tilde{\bm{Q}}_{i,j}(x,y,t; \{\tilde{\bm{V}}^n_{k,l}\}_{\{k,l\} \in {\mathcal{T}}_{C_{i,j}}}).
\end{displaymath}
The mid-point rule is employed again for the quadrature formulas.

In the first-order case, the reconstruction operators are defined as
\begin{subequations}
  \begin{align*}
    \bm{Q}^n_{i,j}(x,y)
     & = \bm{U}^n_{i,j} - \bm{U}^{s,n}_{i,j}(x_i, y_j), \\
    \tilde{\bm{Q}}_{i,j}(x,y,t)
     & = \bm{U}_{i,j}(t) - \bm{U}^n_{i,j}.
  \end{align*}
\end{subequations}
Consequently, the first-order exactly well-balanced reconstruction operator reads
\begin{equation}\label{eq:first-order-wb-reconstruction}
  \bm{P}^{o1}_{i,j}(x,y, t) = \bm{U}^{s,n}_{i,j}(x,y) - \bm{U}^{s,n}_{i,j}(x_i, y_j) + \bm{U}_{i,j}(t),
\end{equation}
while for the second-order scheme, we consider
\begin{subequations}
  \begin{align*}
    \bm{Q}^n_{i,j}(x,y)
     & = \bm{U}^n_{i,j} - \bm{U}^{s,n}_{i,j}(x_i, y_j) + \bm{DV}^n_1 (x-x_i) + \bm{DV}^n_2 (y-y_j),        \\
    \tilde{\bm{Q}}_{i,j}(x,y,t)
     & = \bm{U}_{i,j}(t) - \bm{U}^n_{i,j} + \bm{D\tilde{V}}^n_1(t)(x-x_i) + \bm{D\tilde{V}}^n_2(t)(y-y_j),
  \end{align*}
\end{subequations}
so the resulting exactly well-balanced reconstruction operator in this case is
\begin{equation}\label{eq:second-order-wb-reconstruction}
  \begin{split}
    \bm{P}^{o2}_{i,j}(x, y, t)
     & =  \bm{U}^{s,n}_{i,j}(x,y) -\bm{U}^{s,n}_{i,j}(x_i, y_j) + \bm{U}_{i,j}(t) \\
     & \quad + \bm{DV}^n_1 (x-x_i) + \bm{DV}^n_2 (y-y_j)                          \\
     & \quad + \bm{D\tilde{V}}^n_1(t) (x-x_i) + \bm{D\tilde{V}}^n_2(t) (y-y_j).
  \end{split}
\end{equation}
Using the van Leer limiter, we set
\begin{subequations}
  \begin{align*}
    \bm{DV}^n_1      & = \dfrac{1}{\Delta x} (\phi^n_{i+}(\bm{V}^n_{i,j} - \bm{V}^n_{i-1,j}) + \phi^n_{i-}(\bm{V}^n_{i+1,j} - \bm{V}^n_{i,j})), \\
    \bm{D\bm{V}}^n_2 & = \dfrac{1}{\Delta y} (\phi^n_{j+}(\bm{V}^n_{i,j} - \bm{V}^n_{i,j-1}) + \phi^n_{j-}(\bm{V}^n_{i,j+1} - \bm{V}^n_{i,j})),
  \end{align*}
\end{subequations}
with
\begin{displaymath}
  \phi^n_{i-} = \left\{ \begin{array}{lr}
    \dfrac{|d_{i-}|}{|d_{i-}| + |d_{i+}|} & \text{if} \quad |d_{i-}| + |d_{i+}| > 0, \\[3mm]
    0                                     & \text{otherwise},\end{array}\right.
\end{displaymath}

\begin{displaymath}
  \phi^n_{i+} = \left\{ \begin{array}{lr}
    \dfrac{|d_{i+}|}{|d_{i-}| + |d_{i+}|} & \text{if} \quad |d_{i-}| + |d_{i+}| > 0, \\[3mm]
    0                                     & \text{otherwise},\end{array}\right.
\end{displaymath}

\begin{displaymath}
  \phi^n_{j-} = \left\{ \begin{array}{lr}
    \dfrac{|d_{j-}|}{|d_{j-}| + |d_{j+}|} & \text{if} \quad |d_{j-}| + |d_{j+}| > 0, \\[3mm]
    0                                     & \text{otherwise},\end{array}\right.
\end{displaymath}

\begin{displaymath}
  \phi^n_{j+} = \left\{ \begin{array}{lr}
    \dfrac{|d_{j+}|}{|d_{j-}| + |d_{j+}|} & \text{if} \quad |d_{j-}| + |d_{j+}| > 0, \\[3mm]
    0                                     & \text{otherwise},\end{array}\right.
\end{displaymath}
where
\begin{equation*}
  d_{i-} = \bm{V}^n_{i,j} - \bm{V}^n_{i-1,j}, \  d_{i+} = \bm{V}^n_{i+1,j} - \bm{V}^n_{i,j}, \ d_{j-} = \bm{V}^n_{i,j} - \bm{V}^n_{i,j-1}, \ d_{j+} = \bm{V}^n_{i,j+1} - \bm{V}^n_{i,j},
\end{equation*}
and
\begin{displaymath}
  \bm{V}^n_{k,l} = \bm{U}^n_{k,l} - \bm{U}^{s,n}_{i,j}(x_k, y_l).
\end{displaymath}
For time fluctuations we consider
\begin{subequations}
  \begin{align*}
    \bm{D\tilde{V}}^n_1(t)      & = \dfrac{1}{\Delta x} \Big(\tilde{\phi}_{i+}(\bm{\tilde{V}}^n_{i,j}(t) - \bm{\tilde{V}}^n_{i-1,j}(t)) + \tilde{\phi}_{i-}(\bm{\tilde{V}}^n_{i+1,j}(t) - \bm{\tilde{V}}^n_{i,j}(t))\Big), \\
    \bm{D\bm{\tilde{V}}}^n_2(t) & = \dfrac{1}{\Delta y} \Big(\tilde{\phi}_{j+}(\bm{\tilde{V}}^n_{i,j}(t) - \bm{\tilde{V}}^n_{i,j-1}(t)) + \tilde{\phi}_{j-}(\bm{\tilde{V}}^n_{i,j+1}(t) - \bm{\tilde{V}}^n_{i,j}(t))\Big).
  \end{align*}
\end{subequations}
To prevent the introduction of additional non-linearities, we choose $\tilde{\phi}_{i\pm} = \phi^n_{i\pm}$ and  $\tilde{\phi}_{j\pm} = \phi^n_{j\pm}$.

Finally, in order to approximate the source terms of \eqref{eq:semi-discrete-pressure-system-wb-horizontal} and \eqref{eq:semi-discrete-pressure-system-wb-vertical}, which are time-independent, we use a exactly well-balanced reconstruction operator $P_{i,j,h}$ for the variable $h$ frozen in time $t=t^n$. Then,
\begin{equation}\label{eq:source_terms}
  \begin{array}{rcl}
    ({\mathcal{S}}_1)_{i,j}^{n} & = & \dfrac{1}{|C_{i,j}|}
    \displaystyle\int_{C_{i,j}}
    g(P_{i,j,h}(x,y,t^n) - h^{s,n}_{i,j}(x,y))\,\partial_x H(x,y) \, dy\, dx, \\[5mm]
    ({\mathcal{S}}_2)_{i,j}^{n} & = & \dfrac{1}{|C_{i,j}|}
    \displaystyle\int_{C_{i,j}}
    g(P_{i,j,h}(x,y,t^n) - h^{s,n}_{i,j}(x,y))\,\partial_y H(x,y) \, dy\, dx. \\
  \end{array}
\end{equation}
Furthermore, the interface approximations of the Riemann invariants in systems (\ref{eq:semi-discrete-pressure-system-wb-horizontal}) and (\ref{eq:semi-discrete-pressure-system-wb-vertical}), namely
\begin{equation*}
  (\wRx)_{i+1/2,j}(t) \approx \wRx(x_{i+1/2},y_j,t), \quad (\wLx)_{i+1/2,j}(t) \approx \wLx(x_{i+1/2},y_j,t),
\end{equation*}
and
\begin{equation*}
  (\wRy)_{i,j+1/2}(t) \approx \wRy(x_i,y_{j+1/2},t), \quad (\wLy)_{i,j+1/2}(t) \approx \wLy(x_i,y_{j+1/2},t),
\end{equation*}
correspond to the numerical fluxes at the cell interfaces and are evaluated using a chosen reconstruction operator, whereas the stationary Riemann invariants do not depend on the velocity and are given by
\begin{displaymath}
  (\wRx)^{s,n}_{i,j}(x, y) = (\wLx)^{s,n}_{i,j}(x, y) = (\wRy)^{s,n}_{i,j}(x, y) = (\wLy)^{s,n}_{i,j}(x, y) = \dfrac{g}{2} (h^{s,n}_{i,j}(x, y))^2.
\end{displaymath}

It is possible to recover the pressure and velocities at the interfaces from the Riemann invariants as follows
\begin{equation}\label{pi-star-horizontal}
  \pi^*_{i+1/2,j}(t) = \dfrac{P_{i,j,\wRx}(x_{i+1/2}, y_j, t) + P_{i+1,j,\wLx}(x_{i+1/2}, y_j, t)}{2},
\end{equation}
\begin{equation}\label{u-star}
  u^*_{i+1/2,j}(t) = \dfrac{P_{i,j,\wRx}(x_{i+1/2}, y_j, t) - P_{i+1,j,\wLx}(x_{i+1/2}, y_j, t)}{2a},
\end{equation}
\begin{equation}\label{pi-star-vertical}
  \pi^*_{i,j+1/2}(t) = \dfrac{P_{i,j,\wRy}(x_{i}, y_{j+1/2}, t) + P_{i,j+1,\wLy}(x_{i}, y_{j+1/2}, t)}{2},
\end{equation}
\begin{equation}\label{v-star}
  v^*_{i,j+1/2}(t) = \dfrac{P_{i,j,\wRy}(x_{i}, y_{j+1/2}, t) - P_{i,j+1,\wLy}(x_{i}, y_{j+1/2}, t)}{2a},
\end{equation}
where $P_{i,j,X}$ correspond to specific in-cell reconstruction operators for the variables $X \in \{\wRx, \wLx,\wRy,\wLy \}$.

Conversely, $X^*_{i\pm1/2,j}$ and $X^*_{i,j\pm1/2}$ with $X \in \{h, \qu, \qv, \Pi\}$ in system (\ref{eq:semi-discrete-transport-system}) are reconstructed upwind from (\ref{u-star}) and (\ref{v-star}) as
\begin{displaymath}
  X^*_{i+1/2,j}(t) =
  \left\{ \begin{array}{rcl}
    P_{i,j,X}(x_{i+1/2}, y_j, t)   & \text{ if } & u^*_{i+1/2,j}(t) \geq 0, \\[3mm]
    P_{i+1,j,X}(x_{i+1/2}, y_j, t) & \text{ if } & u^*_{i+1/2,j}(t) < 0,
  \end{array}
  \right.\\
\end{displaymath}
and
\begin{displaymath}
  X^*_{i,j+1/2}(t) =
  \left\{ \begin{array}{rcl}
    P_{i,j,X}(x_{i}, y_{j+1/2}, t)   & \text{ if } & v^*_{i,j+1/2}(t) \geq 0, \\[3mm]
    P_{i,j+1,X}(x_{i}, y_{j+1/2}, t) & \text{ if } & v^*_{i,j+1/2}(t) < 0.
  \end{array}
  \right.
\end{displaymath}
Finally, note that restricting ourselves to schemes of at most second-order allows us to identify cell averages with point values at the cell centres. Consequently, the source terms \eqref{eq:source_terms} vanish, and \eqref{eq:semi-discrete-pressure-system-wb-horizontal} and \eqref{eq:semi-discrete-pressure-system-wb-vertical} lead to the following systems, respectively
\begin{equation}\label{eq:semi-discrete-pressure-system-wb-horizontal-no-source}
  \left\{ \begin{array}{rcl}
    {(\wLx)}^\prime_{i,j}(t) & = & + \dfrac{a}{h_{i,j}(t^n) \Delta x} \Big({(\wLx)}_{i+1/2,j}(t) - (\wLx)^{s,n}_{i,j}(x_{i+1/2}, y_j)\Big) \\[3mm]
                             &   & - \dfrac{a}{h_{i,j}(t^n) \Delta x} \Big({(\wLx)}_{i-1/2,j}(t) - (\wLx)^{s,n}_{i,j}(x_{i-1/2}, y_j)\Big) \\[3mm]
                             &   & -\dfrac{a^2}{h_{i,j}(t^n) \Delta y} \Big({v}^{*}_{i,j+1/2}(t)- {v}^{*}_{i,j-1/2}(t) \Big),              \\[3mm]

    {(\wRx)}^\prime_{i,j}(t) & = & - \dfrac{a}{h_{i,j}(t^n) \Delta x} \Big({(\wRx)}_{i+1/2,j}(t) - (\wRx)^{s,n}_{i,j}(x_{i+1/2}, y_j)\Big) \\[3mm]
                             &   & + \dfrac{a}{h_{i,j}(t^n) \Delta x} \Big({(\wRx)}_{i-1/2,j}(t) - (\wRx)^{s,n}_{i,j}(x_{i-1/2}, y_j)\Big) \\[3mm]
                             &   & - \dfrac{a^2}{h_{i,j}(t^n) \Delta y} \Big({v}^{*}_{i,j+1/2}(t)- {v}^{*}_{i,j-1/2}(t) \Big),             \\[3mm]
  \end{array}\right.
\end{equation}
and
\begin{equation}\label{eq:semi-discrete-pressure-system-wb-vertical-no-source}
  \left\{ \begin{array}{rcl}
    {(\wLy)}^\prime_{i,j}(t) & = & + \dfrac{a}{h_{i,j}(t^n) \Delta y} \Big({(\wLy)}_{i,j+1/2}(t) - (\wLy)^{s,n}_{i,j}(x_{i}, y_{j+1/2})\Big) \\[3mm]
                             &   & - \dfrac{a}{h_{i,j}(t^n) \Delta y} \Big({(\wLy)}_{i,j-1/2}(t) - (\wLy)^{s,n}_{i,j}(x_{i}, y_{j-1/2})\Big) \\[3mm]
                             &   & - \dfrac{a^2}{h_{i,j}(t^n) \Delta x} \Big({u}^{*}_{i+1/2,j}(t)- {u}^{*}_{i-1/2,j}(t)) \Big),              \\[3mm]

    {(\wRy)}^\prime_{i,j}(t) & = & - \dfrac{a}{h_{i,j}(t^n) \Delta y} \Big({(\wRy)}_{i,j+1/2}(t) - (\wRy)^{s,n}_{i,j}(x_{i}, y_{j+1/2})\Big) \\[3mm]
                             &   & + \dfrac{a}{h_{i,j}(t^n) \Delta y} \Big({(\wRy)}_{i,j-1/2}(t) - (\wRy)^{s,n}_{i,j}(x_{i}, y_{j-1/2})\Big) \\[3mm]
                             &   & - \dfrac{a^2}{h_{i,j}(t^n) \Delta x} \Big({u}^{*}_{i+1/2,j}(t)- {u}^{*}_{i-1/2,j}(t)) \Big).              \\[3mm]
  \end{array}\right.
\end{equation}

Having introduced the exactly well-balanced reconstruction operators of first and second order, we now turn to establishing that the resulting semi-discrete schemes are indeed exactly well-balanced. The following proposition shows that, when the numerical solution coincides with a water-at-rest steady state at the beginning of a time step, both the pressure and the transport systems leave it unchanged, for either order of accuracy.

\begin{proposition}
  The schemes that result after considering the semi-discrete schemes \eqref{eq:semi-discrete-pressure-system-wb-horizontal-no-source}, \eqref{eq:semi-discrete-pressure-system-wb-vertical-no-source} and \eqref{eq:semi-discrete-transport-system} with the reconstructions \eqref{eq:first-order-wb-reconstruction} and \eqref{eq:second-order-wb-reconstruction} are exactly well-balanced for water at rest.
\end{proposition}
\textit{Proof}. Let us assume that, at time $t = t^n$, the numerical solution coincides with a discretization of water-at-rest stationary equilibrium $\bm{U}^{s,n}_{i,j}$, that is, a solution verifying the system \eqref{eq:water-at-rest-shallow-water}. Since the pressure is hydrostatic, the system can be equivalently written as
\begin{equation}\label{eq:water-at-rest-eta}
  \left\{ \begin{array}{l}
    u = v = 0,            \\
    h-H = cte.\end{array}\right.
\end{equation}

Let us consider the first-order case. The reconstruction at the interface reads
\begin{equation*}
  \begin{split}
    (\wLx)_{i+1/2,j}(t) & = P^{o1}_{i+1,j,\wLx}(x_{i+1/2},y_j, t)                                                                          \\
                        & =  (\wLx)^{s,n}_{i+1,j}(x_{i+1/2},y_j) - (\wLx)^{s,n}_{i+1,j}(x_{i+1}, y_j) + (\wLx)_{i+1,j}(t)                  \\
                        & =  (\wLx)^{s,n}_{i+1,j}(x_{i+1/2},y_j) - (\wLx)^{s,n}_{i+1,j}(x_{i+1}, y_j) + (\wLx)^{s,n}_{i+1,j}(x_{i+1}, y_j) \\
                        & =  (\wLx)^{s,n}_{i+1,j}(x_{i+1/2},y_j),
  \end{split}
\end{equation*}
where we have used that cell averages coincide with point values at the cell centres, in accordance with \eqref{eq:wb-compatibility-condition} under the mid-point rule. From the second equation of \eqref{eq:water-at-rest-eta} it follows that
\begin{equation*}
  \begin{split}
    h^{s,n}_{i,j}(x_{i+1/2},y_j) & = h^{s,n}_{i+1,j}(x_{i+1/2},y_j), \\
  \end{split}
\end{equation*}
and therefore
\begin{equation*}
  \begin{split}
    (\wLx)^{s,n}_{i+1,j}(x_{i+1/2},y_j) & = \frac{g}{2}\left(h^{s,n}_{i+1,j}(x_{i+1/2},y_j)\right)^2 \\
                                        & = \frac{g}{2}\left(h^{s,n}_{i,j}(x_{i+1/2},y_j)\right)^2   \\
                                        & =  (\wLx)^{s,n}_{i,j}(x_{i+1/2},y_j),
  \end{split}
\end{equation*}
then
\begin{equation*}
  (\wLx)_{i+1/2,j}(t) = (\wLx)^{s,n}_{i,j}(x_{i+1/2},y_j).
\end{equation*}
Let us now turn to the tangential velocities. On the one hand,
\begin{equation*}
  P_{i,j,\wRy}(x_{i}, y_{j+1/2}, t) = (\wRy)^{s,n}_{i,j}(x_i,y_{j+1/2}),
\end{equation*}
while on the other hand
\begin{equation*}
  \begin{split}
    P_{i,j+1,\wLy}(x_{i}, y_{j+1/2}, t) & = (\wLy)^{s,n}_{i,j+1}(x_i, y_{j+1/2}) \\
                                        & = (\wLy)^{s,n}_{i,j}(x_i, y_{j+1/2})   \\
                                        & = (\wRy)^{s,n}_{i,j}(x_i, y_{j+1/2}),  \\
  \end{split}
\end{equation*}
where the second line follows from the third equation of \eqref{eq:water-at-rest-eta} and the third line from the first equation of \eqref{eq:water-at-rest-eta}. Combining both expressions,
\begin{equation*}
  v^*_{i,j+1/2}(t) = \dfrac{P_{i,j,\wRy}(x_{i}, y_{j+1/2}, t) - P_{i,j+1,\wLy}(x_{i}, y_{j+1/2}, t)}{2a}  = 0.
\end{equation*}

Consider now the second-order case. Since the initial condition is stationary, the second and third equations of \eqref{eq:water-at-rest-eta} guarantee that the local steady state $h^{s,n}_{i,j}$ reproduces the data on the stencil direction by direction, that is, $h^{s,n}_{i,j}(x_k,y_l) = h^n_{k,l}$ for every $\{k,l\} \in \mathcal{T}_{C_{i,j}}$; the same holds for each Riemann invariant, as they depend on $h$ only. The spatial fluctuations \eqref{eq:spatial-fluctuation} therefore vanish identically on the stencil, giving $\bm{DV}^n_1 = \bm{DV}^n_2 = 0$. Likewise, evaluating the scheme at the candidate steady state makes the time fluctuations \eqref{eq:time-fluctuation} vanish identically for every $t\in[t^n,t^{n+1}]$ then $\bm{D\tilde{V}}^n_1(t) = \bm{D\tilde{V}}^n_2(t) = 0$. The second-order reconstruction \eqref{eq:second-order-wb-reconstruction} thus reduces to the first-order one \eqref{eq:first-order-wb-reconstruction}, and the computation above applies without modification. Finally, the first equation of system \eqref{eq:semi-discrete-pressure-system-wb-horizontal-no-source} reduces to
\begin{equation*}
  (\wLx)^\prime_{i,j}(t) = 0.
\end{equation*}
The same argument applies to the remaining Riemann invariants and interfaces. Therefore
\begin{equation*}
  \left\{ \begin{array}{rcl}
    {(\wLx)}^\prime_{i,j}(t) & = & 0, \\[3mm]

    {(\wRx)}^\prime_{i,j}(t) & = & 0, \\[3mm]

    {(\wLy)}^\prime_{i,j}(t) & = & 0, \\[3mm]

    {(\wRy)}^\prime_{i,j}(t) & = & 0, \\[3mm]
  \end{array}\right.
\hspace{1cm}
  \left\{ \begin{array}{rcl}
    h^\prime_{i,j}(t)       & = & 0, \\[3mm]

    {(\qu)}^\prime_{i,j}(t) & = & 0, \\[3mm]

    {(\qv)}^\prime_{i,j}(t) & = & 0, \\[3mm]

    \Pi^\prime_{i,j}(t)     & = & 0.
  \end{array}\right.
\end{equation*}
and the stationary solution is preserved. $\square$
\section{Fully-discrete schemes}\label{sec:fully-discrete}
In this section, we present the fully discrete, exactly well-balanced schemes. Two distinct approaches could be considered: a fully explicit formulation—where both the pressure and transport systems are updated explicitly—and a semi-implicit approach, in which the pressure system is integrated implicitly.

The pressure system will be consistently advanced using its associated Riemann invariants rather than directly updating system \eqref{eq:semi-discrete-pressure-system-well-balanced}. Specifically, we either advance the Riemann invariants explicitly or update them implicitly via systems \eqref{eq:semi-discrete-pressure-system-wb-horizontal-no-source}–\eqref{eq:semi-discrete-pressure-system-wb-vertical-no-source}. Once this system is solved, the velocity components are recovered from the invariants and used to reconstruct the directional momenta, thus completing the pressure step.

On the other hand, the transport system is always integrated explicitly. The velocity components are recovered analogously, and the remaining conserved variables are reconstructed using an upwind scheme with respect to these velocities as mentioned before.

Note that the variable $\Pi$ is not updated in either subsystem; its sole purpose is to act as an auxiliary variable to define the Riemann invariants within the relaxed framework for the shallow water equations. Furthermore, since the water height $h$ is frozen during each pressure step, $\pi$ can be naturally identified with the hydrostatic pressure at the beginning of each pressure step.

As previously discussed, in low-Froude-number subcritical regimes ($F_r \ll 1$), the acoustic pressure waves impose a severe restriction on the explicit CFL condition. Because system \eqref{eq:balance-law-pressure} can be decoupled via Riemann invariants, treating the pressure system implicitly bypasses this bottleneck. Given that this efficiency is one of the main contributions of this work, we detail only the semi-implicit formulation; the explicit variant follows analogously.

The time-stepping procedure proceeds as follows. Given a set of cell averages $\bm{U}^n_{i,j}$ at time $t^n$, we first integrate the pressure system over the interval $[t^n, t^{n+1}]$ to obtain the intermediate state $\bm{U}^{n+1-}_{i,j}$. Subsequently, using these intermediate quantities as initial conditions, we advance the transport system \eqref{eq:balance-law-transport} over the same interval $[t^n, t^{n+1}]$, thereby yielding the approximate solution $\bm{U}^{n+1}_{i,j}$ at the new time level $t^{n+1}$.

We formulate the implicit version of the pressure system as
\begin{equation}\label{eq:fully-discrete-pressure-system}
  \left\{ \begin{array}{rcl}
    h^{n+1-}_{i,j}       & = & h^{n}_{i,j},                \\[1mm]
    {(\qu)}^{n+1-}_{i,j} & = & h^{n}_{i,j} u^{n+1-}_{i,j}, \\[1mm]
    {(\qv)}^{n+1-}_{i,j} & = & h^{n}_{i,j} v^{n+1-}_{i,j}, \\
  \end{array}\right.
\end{equation}
where
\begin{displaymath}
  u^{n+1-}_{i,j} = \dfrac{(\wRx)^{n+1-}_{i,j} - (\wLx)^{n+1-}_{i,j}}{2a},
\end{displaymath}
with $(\wRx)^{n+1-}_{i,j}$ and $(\wLx)^{n+1-}_{i,j}$ obtained from system (\ref{eq:semi-discrete-pressure-system-wb-horizontal-no-source}), and
\begin{displaymath}
  v^{n+1-}_{i,j} = \dfrac{(\wRy)^{n+1-}_{i,j} - (\wLy)^{n+1-}_{i,j}}{2a},
\end{displaymath}
with $(\wRy)^{n+1-}_{i,j}$ and $(\wLy)^{n+1-}_{i,j}$ derived from system (\ref{eq:semi-discrete-pressure-system-wb-vertical-no-source}). The fully discrete systems corresponding to (\ref{eq:semi-discrete-pressure-system-wb-horizontal-no-source}) and (\ref{eq:semi-discrete-pressure-system-wb-vertical-no-source}) are given by
\begin{equation}\label{eq:fully-discrete-ri-horizontal-full}
  \left\{ \begin{array}{rcl}
    {(\wLx)}^{n+1-}_{i,j} & = & {(\wLx)}^{n}_{i,j} + \dfrac{a \Delta t}{h^{n}_{i,j} \Delta x} \Big({(\wLx)}^{n+1-}_{i+1/2,j} - (\wLx)^{s,n}_{i,j}(x_{i+1/2}, y_j)\Big) \\[3mm]
                          &   & - \dfrac{a \Delta t}{h^{n}_{i,j} \Delta x} \Big({(\wLx)}^{n+1-}_{i-1/2,j} - (\wLx)^{s,n}_{i,j}(x_{i-1/2}, y_j)\Big)                    \\[3mm]
                          &   & -\dfrac{a^2 \Delta t}{h^{n}_{i,j} \Delta y} \Big({v}^{*,n+1-}_{i,j+1/2}- {v}^{*,n+1-}_{i,j-1/2} \Big),                                 \\[3mm]

    {(\wRx)}^{n+1-}_{i,j} & = & {(\wRx)}^{n}_{i,j} - \dfrac{a \Delta t}{h^{n}_{i,j} \Delta x} \Big({(\wRx)}^{n+1-}_{i+1/2,j} - (\wRx)^{s,n}_{i,j}(x_{i+1/2}, y_j)\Big) \\[3mm]
                          &   & + \dfrac{a \Delta t}{h^{n}_{i,j} \Delta x} \Big({(\wRx)}^{n+1-}_{i-1/2,j} - (\wRx)^{s,n}_{i,j}(x_{i-1/2}, y_j)\Big)                    \\[3mm]
                          &   & - \dfrac{a^2 \Delta t}{h^{n}_{i,j} \Delta y} \Big({v}^{*,n+1-}_{i,j+1/2}- {v}^{*,n+1-}_{i,j-1/2} \Big),                                \\[3mm]
  \end{array}\right.
\end{equation}
and
\begin{equation}\label{eq:fully-discrete-ri-vertical-full}
  \left\{ \begin{array}{rcl}
    {(\wLy)}^{n+1-}_{i,j} & = & {(\wLy)}^{n}_{i,j} + \dfrac{a \Delta t}{h^{n}_{i,j} \Delta y} \Big({(\wLy)}^{n+1-}_{i,j+1/2} - (\wLy)^{s,n}_{i,j}(x_{i}, y_{j+1/2})\Big) \\[3mm]
                          &   & - \dfrac{a \Delta t}{h^{n}_{i,j} \Delta y} \Big({(\wLy)}^{n+1-}_{i,j-1/2} - (\wLy)^{s,n}_{i,j}(x_{i}, y_{j-1/2})\Big)                    \\[3mm]
                          &   & - \dfrac{a^2 \Delta t}{h^{n}_{i,j} \Delta x} \Big({u}^{*,n+1-}_{i+1/2,j}- {u}^{*,n+1-}_{i-1/2,j}) \Big),                                 \\[3mm]

    {(\wRy)}^{n+1-}_{i,j} & = & {(\wRy)}^{n}_{i,j} - \dfrac{a \Delta t}{h^{n}_{i,j} \Delta y} \Big({(\wRy)}^{n+1-}_{i,j+1/2} - (\wRy)^{s,n}_{i,j}(x_{i}, y_{j+1/2})\Big) \\[3mm]
                          &   & + \dfrac{a \Delta t}{h^{n}_{i,j} \Delta y} \Big({(\wRy)}^{n+1-}_{i,j-1/2} - (\wRy)^{s,n}_{i,j}(x_{i}, y_{j-1/2})\Big)                    \\[3mm]
                          &   & - \dfrac{a^2 \Delta t}{h^{n}_{i,j} \Delta x} \Big({u}^{*,n+1-}_{i+1/2,j}- {u}^{*,n+1-}_{i-1/2,j}) \Big),                                 \\[3mm]
  \end{array}\right.
\end{equation}
respectively.

Several numerical tests show that incorporating tangential velocity contributions when reconstructing ${u}^{*,n+1-}_{i\pm1/2,j}$ and ${v}^{*,n+1-}_{i,j\pm1/2}$ from \eqref{u-star} and \eqref{v-star} yields no perceptible gain in accuracy (see Figure \ref{fig:circular_dam_break_O1}), while significantly increasing the computational overhead of the explicit schemes. In the implicit schemes the addition of tangential contributions couples the horizontal and vertical system, thereby losing the favourable properties of the decoupled, independent 1D systems. Therefore, following \cite{fernandez2009modelling}, it is possible to neglect the tangential terms in the orthogonal projection of the Riemann invariants system \eqref{eq:riemann-invariants-system} at the cell interfaces from the outset, which yields the following fully discrete schemes:
\begin{equation}\label{eq:fully-discrete-ri-horizontal}
  \left\{ \begin{array}{rcl}
    {(\wLx)}^{n+1-}_{i,j} & = & {(\wLx)}^{n}_{i,j} + \dfrac{a \Delta t}{h^{n}_{i,j} \Delta x} \Big({(\wLx)}^{n+1-}_{i+1/2,j} - (\wLx)^{s,n}_{i,j}(x_{i+1/2}, y_j)\Big) \\[3mm]
                          &   & - \dfrac{a \Delta t}{h^{n}_{i,j} \Delta x} \Big({(\wLx)}^{n+1-}_{i-1/2,j} - (\wLx)^{s,n}_{i,j}(x_{i-1/2}, y_j)\Big),                   \\[3mm]

    {(\wRx)}^{n+1-}_{i,j} & = & {(\wRx)}^{n}_{i,j} - \dfrac{a \Delta t}{h^{n}_{i,j} \Delta x} \Big({(\wRx)}^{n+1-}_{i+1/2,j} - (\wRx)^{s,n}_{i,j}(x_{i+1/2}, y_j)\Big) \\[3mm]
                          &   & + \dfrac{a \Delta t}{h^{n}_{i,j} \Delta x} \Big({(\wRx)}^{n+1-}_{i-1/2,j} - (\wRx)^{s,n}_{i,j}(x_{i-1/2}, y_j)\Big),                   \\[3mm]
  \end{array}\right.
\end{equation}
\begin{equation}\label{eq:fully-discrete-ri-vertical}
  \left\{ \begin{array}{rcl}
    {(\wLy)}^{n+1-}_{i,j} & = & {(\wLy)}^{n}_{i,j} + \dfrac{a \Delta t}{h^{n}_{i,j} \Delta y} \Big({(\wLy)}^{n+1-}_{i,j+1/2} - (\wLy)^{s,n}_{i,j}(x_{i}, y_{j+1/2})\Big) \\[3mm]
                          &   & - \dfrac{a \Delta t}{h^{n}_{i,j} \Delta y} \Big({(\wLy)}^{n+1-}_{i,j-1/2} - (\wLy)^{s,n}_{i,j}(x_{i}, y_{j-1/2})\Big),                   \\[3mm]

    {(\wRy)}^{n+1-}_{i,j} & = & {(\wRy)}^{n}_{i,j} - \dfrac{a \Delta t}{h^{n}_{i,j} \Delta y} \Big({(\wRy)}^{n+1-}_{i,j+1/2} - (\wRy)^{s,n}_{i,j}(x_{i}, y_{j+1/2})\Big) \\[3mm]
                          &   & + \dfrac{a \Delta t}{h^{n}_{i,j} \Delta y} \Big({(\wRy)}^{n+1-}_{i,j-1/2} - (\wRy)^{s,n}_{i,j}(x_{i}, y_{j-1/2})\Big).                   \\[3mm]
  \end{array}\right.
\end{equation}
Note that the interface values ${(\wLx)}^{n+1-}_{i+1/2,j}$ and ${(\wLy)}^{n+1-}_{i,j+1/2}$ are evaluated as
\begin{displaymath}
  {(\wLx)}^{n+1-}_{i+1/2,j} = P_{i+1,j,\wLx}(x_{i+1/2}, y_j, t^{n+1-}), {(\wLy)}^{n+1-}_{i,j+1/2} = P_{i,j+1,\wLy}(x_{i}, y_{j+1/2}, t^{n+1-}),
\end{displaymath}
once $\wLx$ and $\wLy$ have been computed at the cell centres. All remaining terms are defined analogously.

Consequently, the implicit pressure step reduces to solving a set of decoupled linear systems, one for each of the Riemann invariants $\wRx$, $\wLx$, $\wRy$ and $\wLy$, which can be solved very efficiently. The resulting solutions are subsequently used to update \eqref{eq:fully-discrete-pressure-system}. In the present implementation, these sparse linear systems are solved using a direct method.

Because the tangential terms are neglected, the pressure step decouples into
independent one-dimensional problems along each grid line: $2N_y$ systems of
size $N_x$ and $2N_x$ systems of size $N_y$. Each system is bidiagonal for
the first-order reconstruction and has four diagonals for the second-order
one. Since the bandwidth is fixed regardless of the mesh size, each line can
be solved in linear time by a banded LU factorisation, so that the overall
cost of the pressure step is $O(N_xN_y)$, the same order as an explicit
step. This structure has two practical consequences.

First, the cost of the implicit step is independent of the stability
restriction on the pressure system: any increase in the time step directly
reduces the total number of time steps required.

Second, since the systems are mutually independent, the pressure step is
naturally parallelisable. The only additional requirement is that, for the
second-order reconstruction, partial pivoting must be used to preserve
stability.

Periodic boundary conditions introduce only a few entries outside the band,
which are treated as a low-rank correction via the Sherman--Morrison--Woodbury
formula.

The transport system (\ref{eq:semi-discrete-transport-system}) is advanced using $(h^{n+1-}_{i,j}, {(\qu)}^{n+1-}_{i,j}, {(\qv)}^{n+1-}_{i,j})$ as the initial condition. Specifically
\begin{equation}\label{eq:fully-discrete-transport-system}
  \left\{ \begin{array}{rcl}
    h^{n+1}_{i,j}       & = & h^{n+1-}_{i,j} - \dfrac{\Delta t}{\Delta x} \Big(h^{*,n+1-}_{i+1/2,j}{u}^{*,n+1-}_{i+1/2,j}-h^{*,n+1-}_{i-1/2,j}{u}^{*,n+1-}_{i-1/2,j} \Big)                   \\[3mm]
                        &   & - \dfrac{\Delta t}{\Delta y} \Big(h^{*,n+1-}_{i,j+1/2}v^{*,n+1-}_{i,j+1/2} - h^{*,n+1-}_{i,j-1/2}v^{*,n+1-}_{i,j-1/2}\Big),                                    \\[3mm]

    {(\qu)}^{n+1}_{i,j} & = & {(\qu)}^{n+1-}_{i,j} - \dfrac{\Delta t}{\Delta x} \Big({(\qu)}^{*,n+1-}_{i+1/2,j}{u}^{*,n+1-}_{i+1/2,j}-{(\qu)}^{*,n+1-}_{i-1/2,j}{u}^{*,n+1-}_{i-1/2,j} \Big) \\[3mm]
                        &   & -\dfrac{\Delta t}{\Delta y} \Big({(\qu)}^{*,n+1-}_{i,j+1/2}{v}^{*,n+1-}_{i,j+1/2}-{(\qu)}^{*,n+1-}_{i,j-1/2}{v}^{*,n+1-}_{i,j-1/2} \Big),                      \\[3mm]

    {(\qv)}^{n+1}_{i,j} & = & {(\qv)}^{n+1-}_{i,j} - \dfrac{\Delta t}{\Delta x} \Big({(\qv)}^{*,n+1-}_{i+1/2,j}{u}^{*,n+1-}_{i+1/2,j}-{(\qv)}^{*,n+1-}_{i-1/2,j}{u}^{*,n+1-}_{i-1/2,j} \Big) \\[3mm]
                        &   & -\dfrac{\Delta t}{\Delta y} \Big({(\qv)}^{*,n+1-}_{i,j+1/2}{v}^{*,n+1-}_{i,j+1/2}-{(\qv)}^{*,n+1-}_{i,j-1/2}{v}^{*,n+1-}_{i,j-1/2} \Big).                      \\[3mm]
  \end{array}\right.
\end{equation}

Here, the normal interface velocities ${u}^{*,n+1-}_{i+1/2,j}$ and ${v}^{*,n+1-}_{i,j+1/2}$ are computed at each cell edge as follows
\begin{displaymath}
  {u}^{*,n+1-}_{i+1/2,j} = \dfrac{P_{i,j,\wRx}(x_{i+1/2}, y_j, t^{n+1-}) - P_{i+1,j,\wLx}(x_{i+1/2}, y_j, t^{n+1-})}{2a},
\end{displaymath}
where the invariant states $\wRx$ and $\wLx$ are obtained from \eqref{eq:fully-discrete-ri-horizontal-full} or \eqref{eq:fully-discrete-ri-horizontal}, and
\begin{displaymath}
  {v}^{*,n+1-}_{i,j+1/2} = \dfrac{P_{i,j,\wRy}(x_{i}, y_{j+1/2}, t^{n+1-}) - P_{i,j+1,\wLy}(x_{i}, y_{j+1/2}, t^{n+1-})}{2a},
\end{displaymath}
where the corresponding states for $\wRy$ and $\wLy$ are derived from (\ref{eq:fully-discrete-ri-vertical-full}) or (\ref{eq:fully-discrete-ri-vertical}).

Finally, the variables $h^{*,n+1-}_{i+1/2,j}$, $h^{*,n+1-}_{i,j+1/2}$, $(\qu)^{*,n+1-}_{i+1/2,j}$, $(\qu)^{*,n+1-}_{i,j+1/2}$, $(\qv)^{*,n+1-}_{i+1/2,j}$ and $(\qv)^{*,n+1-}_{i,j+1/2}$ are also evaluated using the exactly well-balanced reconstruction operator and the upwind scheme, yielding
\begin{displaymath}
  X^{*,n+1-}_{i+1/2,j} = \left\{ \begin{array}{rcl}
    P_{i,j,X}(x_{i+1/2},y_j,t^{n+1-})   & \text{ if } & {u}^{*,n+1-}_{i+1/2,j} > 0,    \\
    P_{i+1,j,X}(x_{i+1/2},y_j,t^{n+1-}) & \text{ if } & {u}^{*,n+1-}_{i+1/2,j} \leq 0,
  \end{array}\right.
\end{displaymath}
and
\begin{displaymath}
  X^{*,n+1-}_{i,j+1/2} = \left\{ \begin{array}{rcl}
    P_{i,j,X}(x_{i},y_{j+1/2},t^{n+1-})   & \text{ if } & {v}^{*,n+1-}_{i,j+1/2} > 0,    \\
    P_{i,j+1,X}(x_{i},y_{j+1/2},t^{n+1-}) & \text{ if } & {v}^{*,n+1-}_{i,j+1/2} \leq 0,
  \end{array}\right.
\end{displaymath}
where $X \in \{h, \qu, \qv\}$ and the rest are defined analogously.

We employ an operator-splitting technique whereby the pressure system is integrated prior to the transport step, both utilizing a first-order Euler time discretization. While this specific sequence is not strictly mandatory, concluding the advancement with the transport step has been empirically shown to yield better results in previous studies \cite{caballero2024semi}. With the aim of extending the scheme to second order, we
apply a palindromic Strang splitting methodology \cite{marchuk1988metody,marchuk1990splitting,strang1968construction}, structured as follows:
\begin{enumerate}
  \item Perform a transport step over a time step $\frac{\Delta t}{2}$.
  \item Perform a pressure step over a full time step $\Delta t$.
  \item Perform a final transport step over a time step $\frac{\Delta t}{2}$.
\end{enumerate}
Although the intra-step time integrations remain inherently first-order, the second-order spatial reconstructions applied within each sub-step, together with the palindromic structure of the splitting, are the main ingredients through which we seek to achieve second-order accuracy.

\section{Numerical experiments}\label{sec:numerical-experiments}

In this section, we present several numerical experiments to evaluate the accuracy, robustness, and computational efficiency of the proposed schemes. Throughout this section, EXP denotes the fully explicit schemes, whereas IMP refers to the semi-implicit variants, where the pressure is solved implicitly. O1 and O2 denote first- and second-order spatial reconstructions, respectively. Moreover, the acronyms PT, TP, PTP, and TPT specify the sequential order in which the pressure (P) and transport (T) sub-steps are executed.

The stability restriction is determined separately for the pressure and transport sub-steps via the respective Courant–Friedrichs–Lewy parameters, $\text{CFL}_p$ and $\text{CFL}_t$. These restrictions are directly governed by the maximum eigenvalues of the continuous pressure system \eqref{eq:eigenvalues_pressure} and transport system \eqref{eq:eigenvalues_transport} respectively. For the fully explicit schemes (EXP), numerical stability requires both conditions $\text{CFL}_p < 1$ and $\text{CFL}_t < 1$ to be satisfied as usual in two dimensions. In contrast, for the semi-implicit formulations (IMP), treating the pressure step implicitly removes the restrictive acoustic CFL condition associated with gravity waves. Consequently, $\text{CFL}_p$ can be chosen significantly larger than $1$ ($\text{CFL}_p > 1$), yielding substantial efficiency gains while leaving the transport step constrained solely by the physical fluid velocity and $\text{CFL}_t < 1$.

Accordingly, the global time step $\Delta t$ is selected as
\begin{displaymath}
  \Delta t = \min \{ \Delta t_p, \Delta t_t \},
\end{displaymath}
where
\begin{displaymath}
  \Delta t_p = \dfrac{\text{CFL}_p \min \{\Delta x , \Delta y \} }{2\max_{i,j} \left \{\frac{a}{h^n_{i,j}} \right \}}, \quad \Delta t_t = \dfrac{\text{CFL}_t \min \{\Delta x , \Delta y \} }{2\max_{i,j} \left \{\left |u^n_{i,j} \right |, \left |v^n_{i,j} \right | \right \}},
\end{displaymath}
with
\begin{equation*}
  u^n_{i,j} = \frac{(q_u)^n_{i,j}}{h^n_{i,j}}, \quad v_{i,j} = \frac{(q_v)^n_{i,j}}{h^n_{i,j}}
\end{equation*}
and
\begin{displaymath}
  a = \max_{i,j} \left \{ h^n_{i,j} \sqrt{g h^n_{i,j}} \right \}.
\end{displaymath}
In all numerical tests, the gravitational acceleration is set to $g = 9.81$, and the reference water level is taken as zero. Unless specified otherwise, periodic boundary conditions are imposed across all domain boundaries.

\subsection{Circular dam break over a non-flat topography}
The circular dam-break configuration is a standard benchmark for assessing how well a numerical method preserves radial symmetry \cite{alcrudo1993high}. We consider the square domain $\Omega =  [0,50]\times[0,50]$, in which two regions are separated by a cylindrical wall of radius $11$ centred at $(25,25)$.

This test was developed in \cite{aureli2008weighted} to investigate the hydrodynamic response triggered by the sudden collapse of an idealized circular dam over a non-flat bottom. The expected behavior is that when the water starts to drain, the circular dam-break bore waves spread and propagate radially and symmetrically.

\begin{figure}[htbp]
  \centering
  \includegraphics[width=0.7\linewidth]{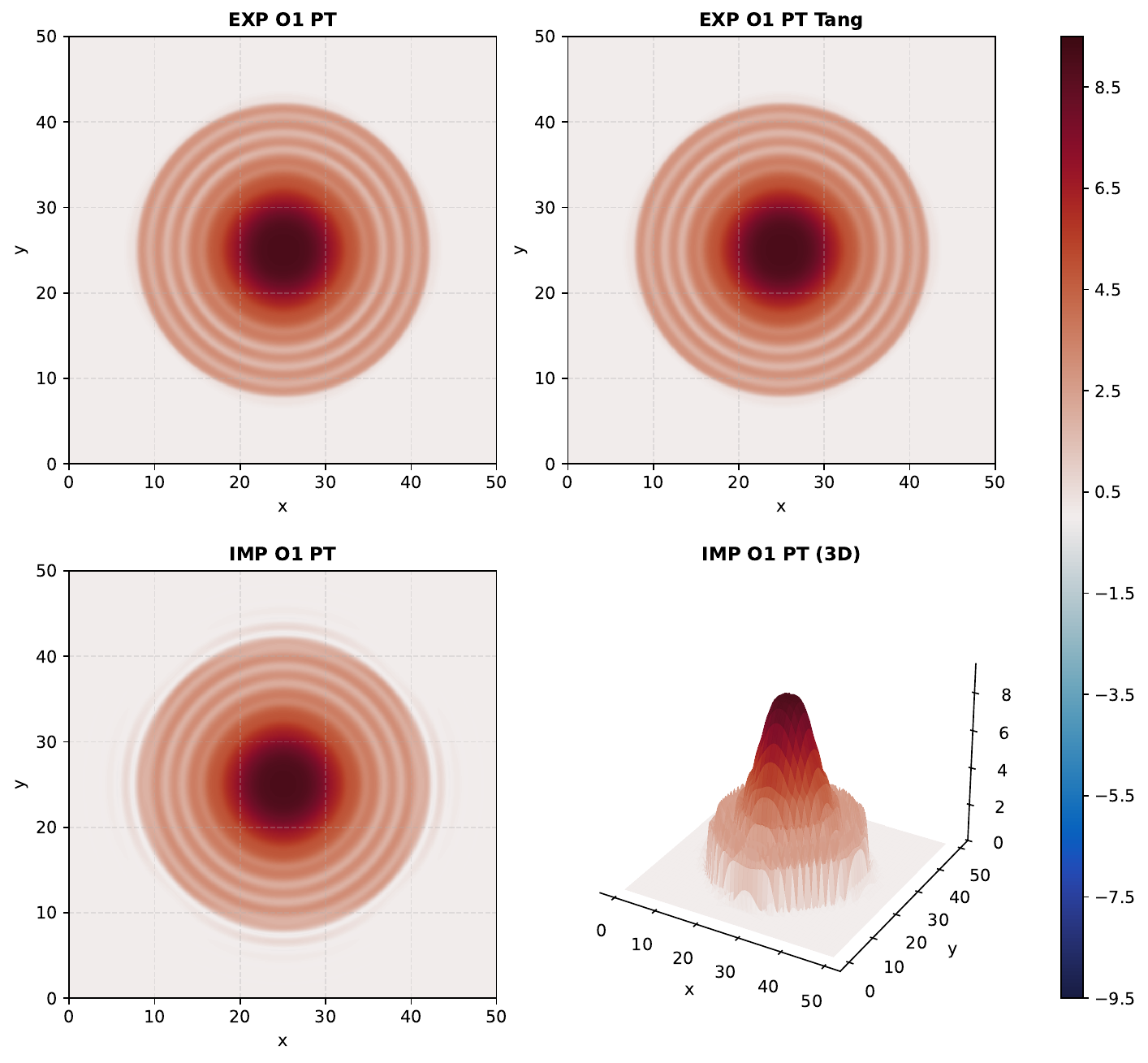}
  \caption{Free surface elevation ($\eta = h - H$) at time $t = 0.69$ for the radial dam-break problem over a non-flat topography ($501 \times 501$ uniform grid). Numerical solutions obtained using two-step first-order (O1) schemes are compared: fully explicit ($\text{CFL}_p = \text{CFL}_t = 1$) with and without tangential velocity contributions (suffix ``Tang'') versus semi-implicit ($\text{CFL}_p = 5$, $\text{CFL}_t = 1$).}
  \label{fig:circular_dam_break_O1}
\end{figure}

\begin{figure}[htbp]
  \centering
  \includegraphics[width=0.7\linewidth]{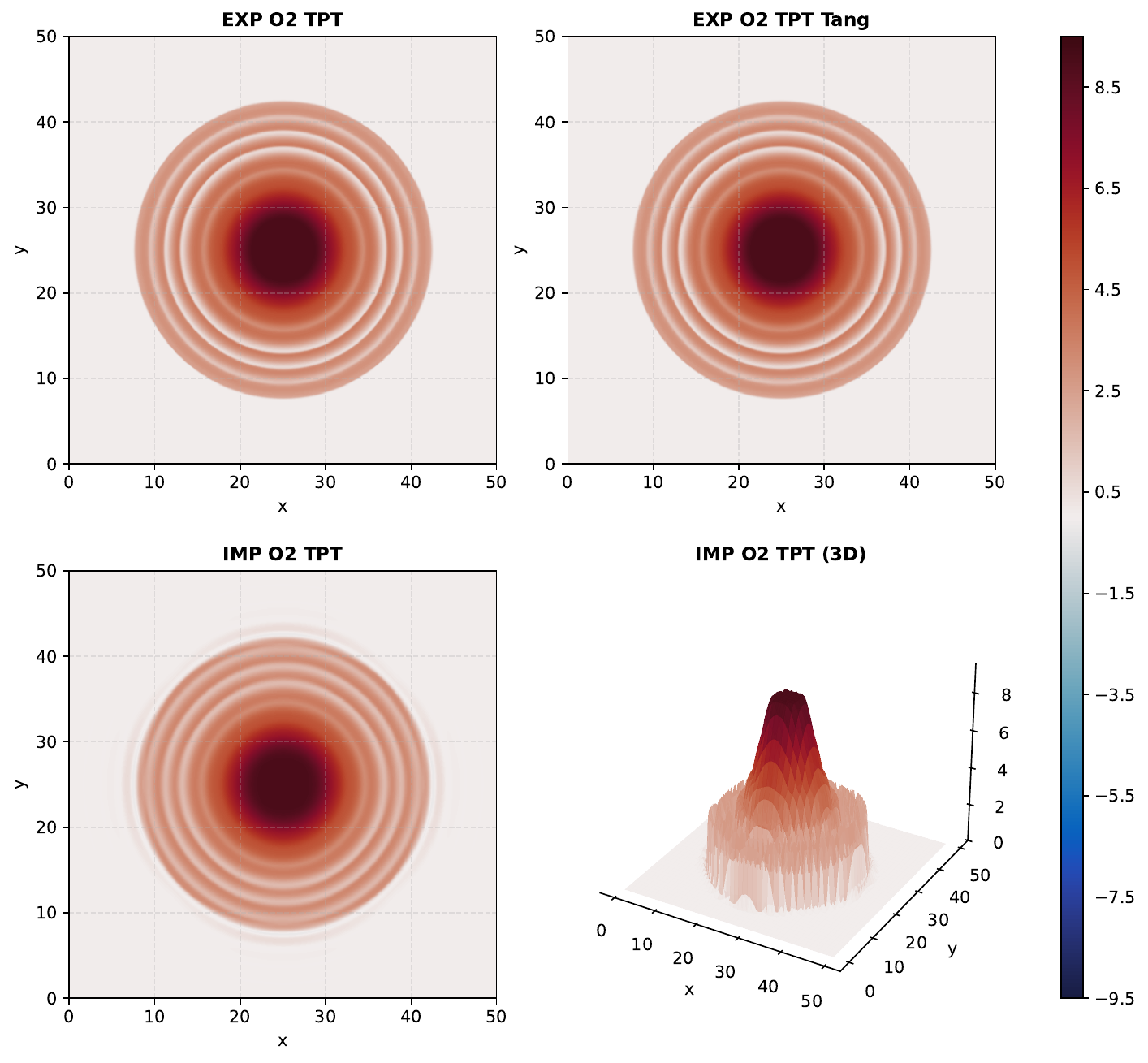}
  \caption{Free surface elevation ($\eta = h - H$) at time $t = 0.69$ for the radial dam-break problem over a non-flat topography ($501 \times 501$ uniform grid). Numerical solutions obtained using three-step second-order (O2) schemes are compared: fully explicit ($\text{CFL}_p = \text{CFL}_t = 1$) including variant with tangential contributions (``Tang'') versus semi-implicit ($\text{CFL}_p = 5$, $\text{CFL}_t = 1$).}
  \label{fig:circular_dam_break_O2}
\end{figure}

The bottom depth is prescribed as
\begin{displaymath}
  H(x,y) = 1 + \frac{1}{2} \Big(1 + \cos \Big(\pi \ \sqrt{(x-25)^2 + (y-25)^2} \Big) \Big)
\end{displaymath}
while the initial velocity field is set to $u \equiv v \equiv 0$ and the initial free surface is prescribed by
\begin{displaymath}
  \eta(x,y,0) = h(x,y,0) - H(x,y) = \begin{cases}
    9 & \text{if } \sqrt{(x-25)^2 + (y-25)^2} \leq 11, \\
    0 & \text{otherwise}.
  \end{cases}
\end{displaymath}
The numerical solution is displayed at $t=0.69$. First-order results are reported in Figure \ref{fig:circular_dam_break_O1}, whereas their second-order counterparts are illustrated in Figure \ref{fig:circular_dam_break_O2}. In both figures, simulations incorporating tangential velocity contributions for the explicit schemes are distinguished by the suffix ``Tang'' in the legend.

All simulations are performed on a uniform $501\times501$ mesh. The enhanced stability of the semi-implicit treatment allows us to bypass the restrictive acoustic Courant limit, safely employing $\text{CFL}_p = 5$ while keeping $\text{CFL}_t = 1$ for the explicit transport stage.

The results of this two-dimensional radial benchmark confirm several key properties of the proposed numerical framework. First, despite the underlying Cartesian mesh discretisation, Figures \ref{fig:circular_dam_break_O1} and \ref{fig:circular_dam_break_O2} show that, for both the first- and second-order schemes, the expanding bore and the waves propagate radially without any grid-orientation bias. Second, the well-balanced treatment robustly handles the interaction between the cylindrical shock wave and the non-trivial bottom topography, preventing spurious oscillations induced by the bottom topography. Finally, and most importantly from a computational efficiency standpoint, the semi-implicit formulation successfully operates at an acoustic Courant number of $\text{CFL}_p = 5$, a fivefold increase over the explicit limit.

\subsection{Square dam break}
\label{test:square-dam-break}

Next, a two-dimensional square dam-break test case is presented to evaluate the capability of the numerical schemes to capture accurately complex hydrodynamic patterns within the velocity field. The computational domain is defined as $\Omega = [-1,1]\times[-1,1]$, which is discretized into a uniform Cartesian grid of $400\times400$ cells. For the explicit schemes, the stability parameters are set to $\text{CFL}_t = \text{CFL}_p = 1$, whereas for the semi-implicit formulations, we assign $\text{CFL}_t = 1$ and a significantly relaxed acoustic Courant number of $\text{CFL}_p = 4$. Initially, the fluid is at rest throughout the domain, i.e., $u \equiv v \equiv 0$.

The bottom topography is set to be flat, $H(x,y) = 1$, and the initial free surface elevation is defined as
\begin{equation*}
\eta(x,y,0) = \begin{cases}
    1 & \text{if } |x| \leq 0.05 \text{ and } |y| \leq 0.05, \\
    0 & \text{otherwise}.
  \end{cases}
\end{equation*}

\begin{figure}[htbp]
  \centering
  \includegraphics[width=0.99\linewidth]{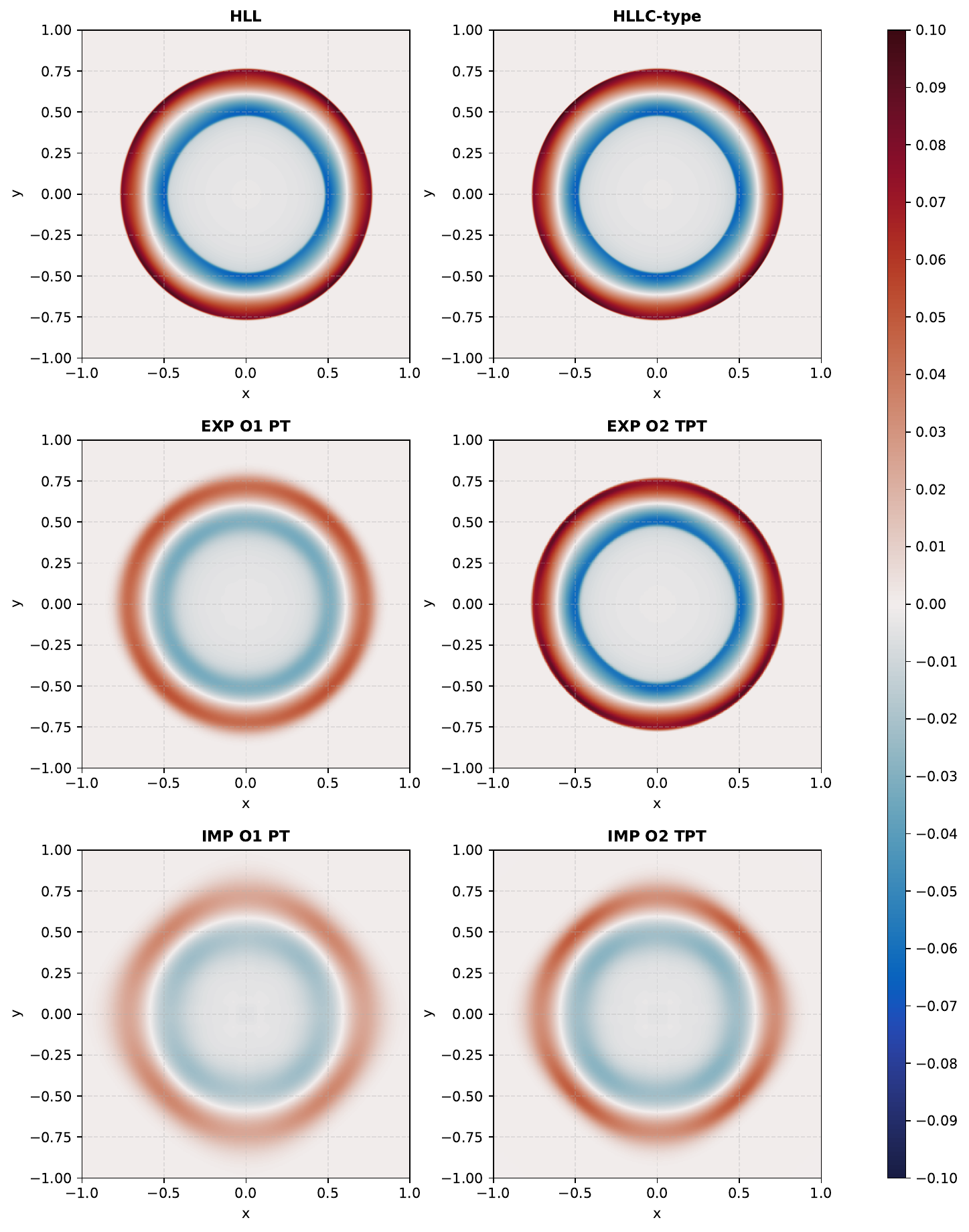}
  \caption{Square dam-break problem: Comparison of the free surface elevation ($\eta = h - H$) at time $t=0.2$. The simulations of the proposed explicit and semi-implicit schemes are performed on a $400\times400$ grid, whereas the HLL and HLLC-type reference solutions are computed on a fine $3200\times3200$ mesh to establish a baseline. All schemes correctly capture the macro-scale outward-propagating shock waves and rarefaction fans.}
  \label{fig:sdb-eta-comparison}
\end{figure}

\begin{figure}[htbp]
  \centering
  \includegraphics[width=1.0\linewidth]{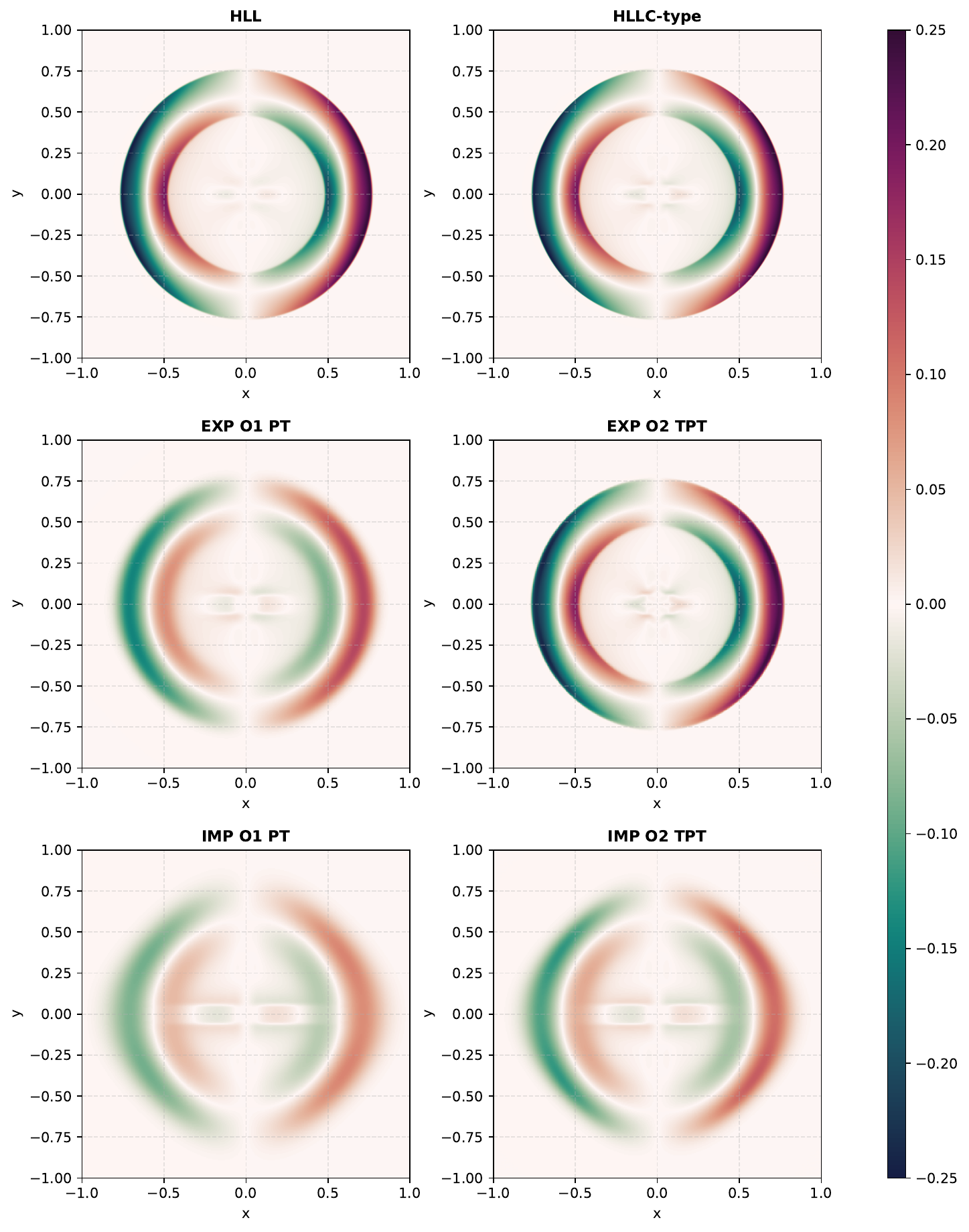}
  \caption{Square dam-break problem: Horizontal velocity field ($u=\frac{q_u}{h}$) at time $t=0.2$. Our proposed explicit and semi-implicit schemes ($400\times400$ grid) are compared against fine-mesh ($3200\times3200$) HLL and HLLC-type solutions. Notably, the standard HLL scheme requires a highly refined grid to display the central structure due to numerical smearing on coarse grids. In contrast, the proposed schemes capture the central structure on the coarser mesh, matching the HLLC-type reference.}
  \label{fig:sdb-velocity-u}
\end{figure}

\begin{figure}[htbp]
  \centering
  \includegraphics[width=1.0\linewidth]{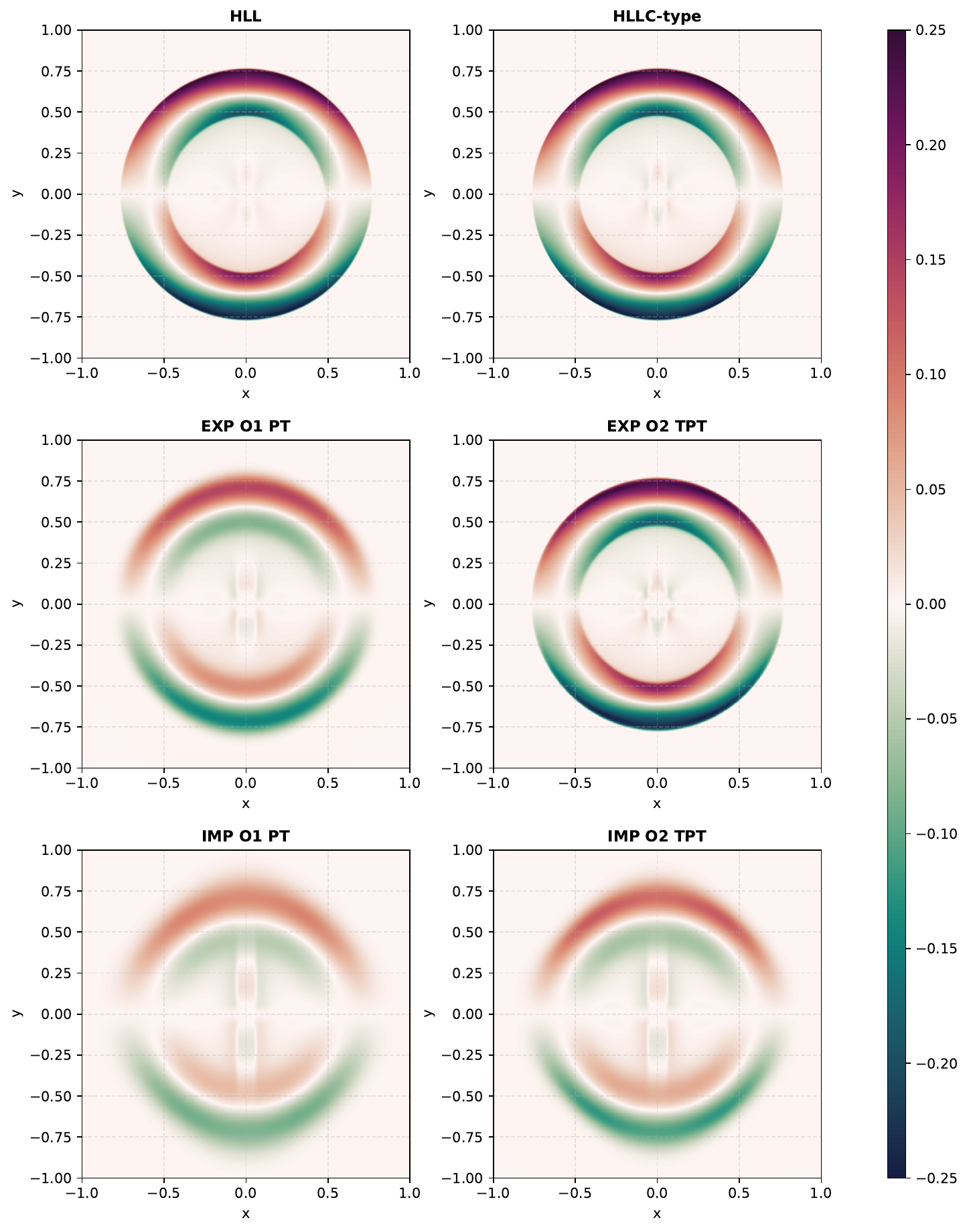}
  \caption{Square dam-break problem: Vertical velocity field ($v = \frac{q_v}{h}$) at time $t=0.2$. Comparison between the proposed schemes on a $400\times400$ grid and reference HLL and HLLC-type solvers on a refined $3200\times3200$ mesh. The proposed semi-implicit methodology demonstrates the capacity to resolve complex inner hydrodynamic structures without smearing, closely mirroring the contact-preserving HLLC-type behavior without requiring excessive grid refinement.}
  \label{fig:sdb-velocity-v}
\end{figure}

While the overall evolution of the free surface elevation behaves as theoretically expected (see Figure \ref{fig:sdb-eta-comparison}), a distinct structure is observed in the velocity field at the centre of the domain (see Figures \ref{fig:sdb-velocity-u} and \ref{fig:sdb-velocity-v}).

\begin{figure}[htbp]
  \centering
  \includegraphics[width=\linewidth]{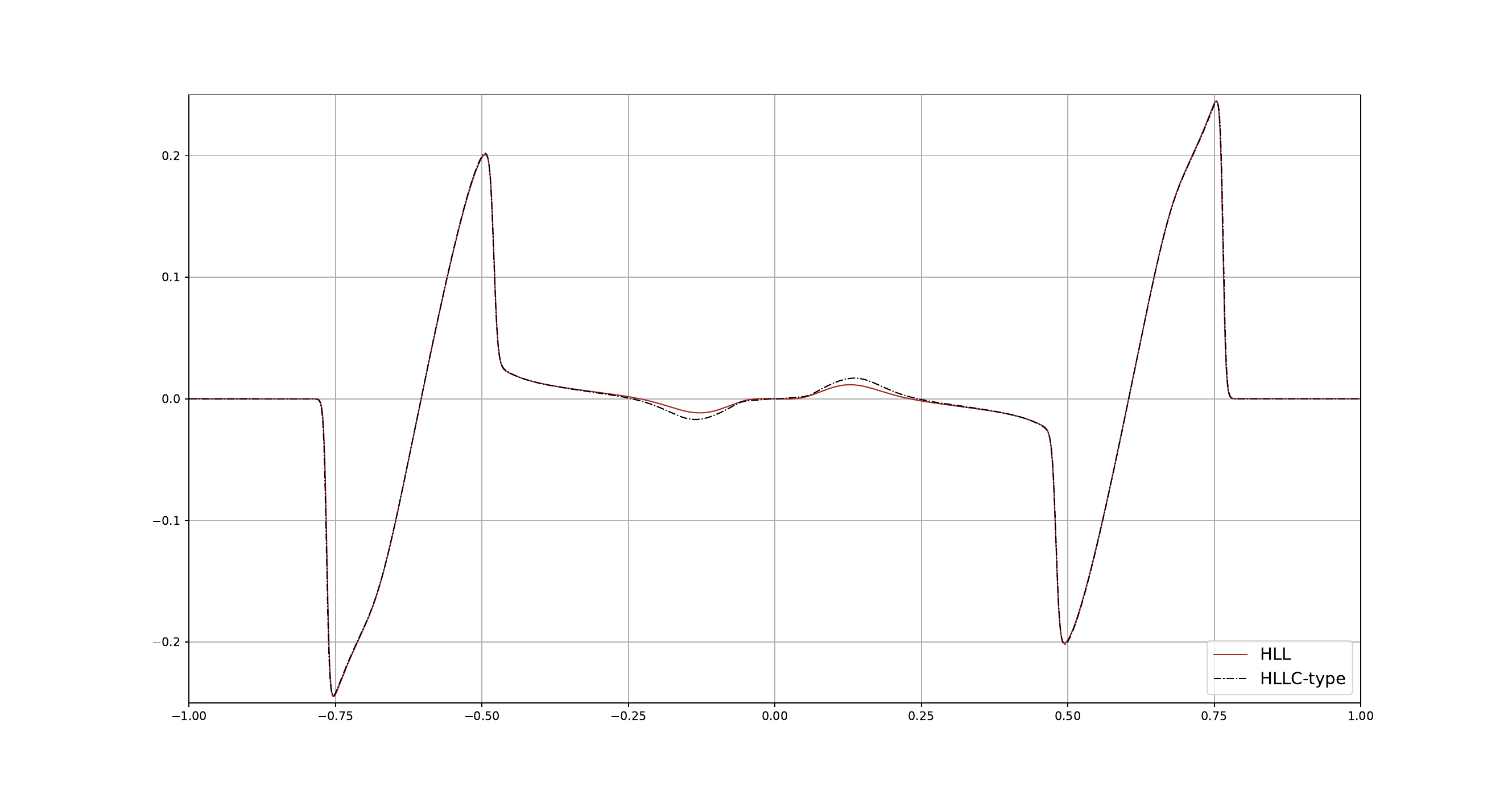}
  \caption{Square dam-break problem: horizontal cut ($y \approx 0.0022$) of the
    horizontal velocity field ($u = \frac{q_u}{h}$) at time $t=0.2$.
    Comparison between the reference solutions ($3200 \times 3200$) obtained
    with the HLL and HLLC-type solvers.}
  \label{fig:sdb-u-cut-refrences}
\end{figure}

\begin{figure}[htbp]
  \centering
  \includegraphics[width=\linewidth]{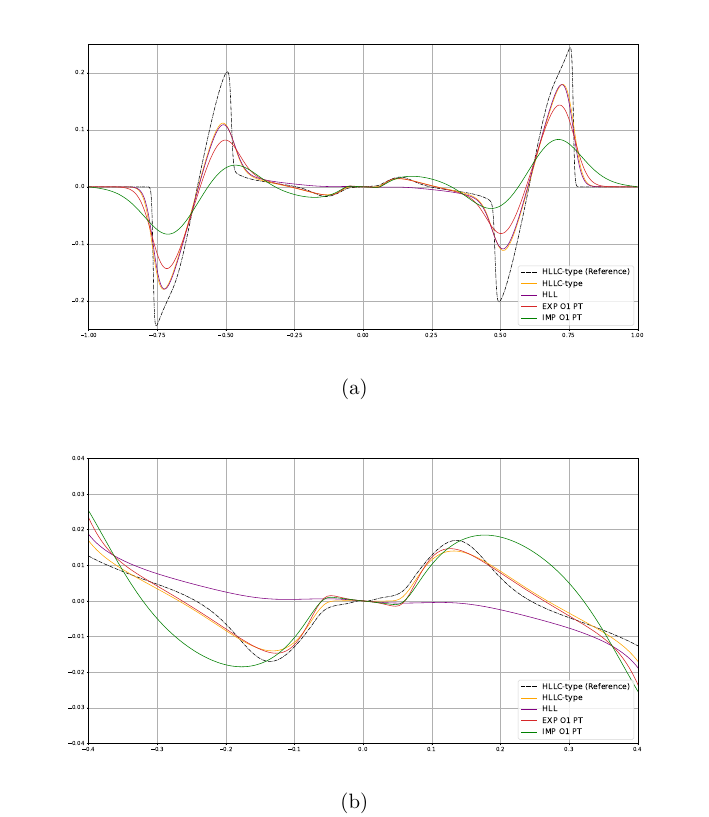}
  \caption{Square dam-break problem: horizontal cut ($y \approx 0.0025$) comparing
    the first-order schemes at coarse resolution ($400 \times 400$) with the
    HLLC-type solver from Figure \ref{fig:sdb-u-cut-refrences}. (a)
    shows the whole domain, while (b) shows a zoom of the
    central structure.}
  \label{fig:sdb-u-cut-comparison}
\end{figure}

Far from being a numerical artifact, this central structure is accurately recovered by standard first-order HLL \cite{toro2009riemann} solvers only when the grid is sufficiently refined, as well as by a first-order HLLC-type scheme \cite{de2013efficient,fernandez2008consistent} (see Figure \ref{fig:sdb-u-cut-refrences}), which inherently provides a superior resolution of contact discontinuities. To properly contextualize our results, the proposed schemes are simulated on a $400\times400$ grid and compared against reference HLL and HLLC-type solutions computed on a highly refined $3200\times3200$ mesh. While an HLLC-type solver is capable of resolving this central hydrodynamic structure even on coarse grids, the standard HLL Riemann solver suffers from excessive numerical smearing at coarse resolutions (see Figure \ref{fig:sdb-u-cut-comparison}), failing to capture the inner feature unless a very fine grid is used. Remarkably, our proposed semi-implicit schemes—despite operating on the $400\times400$ mesh—preserve this central flow pattern with clarity and sharpness, achieving a structural fidelity that closely mirrors the high-resolution HLLC-type reference solution.

In summary, this test shows that the central structure in the velocity field is a genuine feature of the solution: it is recovered by two independent solvers on a fine mesh and persists under refinement. The proposed schemes resolve it on a $400\times400$ mesh with an accuracy that the HLL solver only attains on a finer mesh, and comparable to that of the HLLC-type solver. This is consistent with the structure of the method: the acoustic sub-step does not act on the tangential velocity, which is only transported upwind by the interface normal velocity, so that no additional diffusion is introduced on transverse velocity gradients.

\subsection{Flow over three mounds}

\begin{figure}[htbp]
  \centering
  \includegraphics[width=0.96\linewidth]{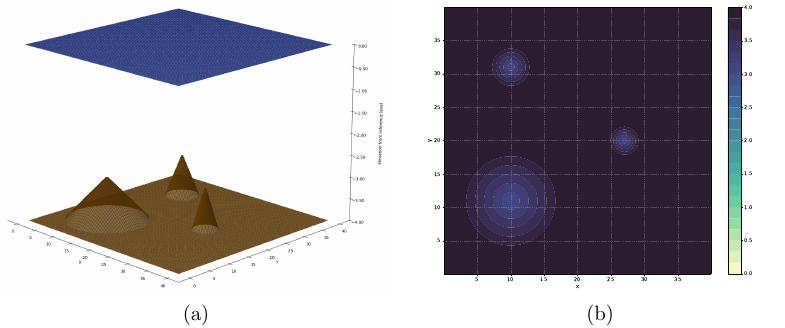}
  \caption{Initial water at rest equilibrium state over a non-trivial three-mound bathymetry on the domain $\Omega = [0,40] \times [0,40]$. (a) Free surface elevation ($\eta = h - H = 0$) and topography ($z=-H$). (b) Corresponding bottom depth ($H$) measured from the reference mean sea level.}
  \label{fig:well_balanced_initial_condition}
\end{figure}

The next experiment is designed to evaluate the well-balanced character of the proposed numerical schemes under constraints imposed by a non-trivial bottom topography. On the domain $\Omega =  [0,40]\times[0,40]$, the bottom topography includes three mounds. The bottom depth, measured from the reference level, is given by
\begin{equation*}
  \begin{split}
    H(x,y) = -\max\Big(-4, \max \Big(-3- \frac{1}{8}\sqrt{(x-10)^2 + (y-11)^2}, \\
    \max \Big(-3-\frac{3}{10}\sqrt{(x-10)^2 + (y-31)^2},                        \\
        -3-\frac{4}{10}\sqrt{(x-27)^2 + (y-20)^2} \Big) \Big) \Big),
  \end{split}
\end{equation*}
and is depicted in Figure \ref{fig:well_balanced_initial_condition}. The system is initialized in a water at rest equilibrium state, setting the free surface elevation and velocity fields to zero ($\eta \equiv u \equiv v \equiv 0$). The simulation is carried out up to the final time $T=200$ on a uniform $200\times200$ grid.
As reported in Table \ref{nr:wb-exp}, the explicit schemes preserve the water at rest equilibrium exactly ($0.00\text{e}+00$) across all variables, confirming that the exactly well-balanced spatial reconstruction perfectly balances the divergence of the fluxes and the source terms at the discrete level. Notably, the semi-implicit formulations maintain the equilibrium state with negligible $L^1$-norm relative errors \footnote{Defined as
$\frac{\|q(x,y,T) - q(x,y,0)\|_{L^1}}{|\Omega|}$, where $q$ denotes the
conserved quantity under consideration.} of order $10^{-13}$ to $10^{-15}$ (Table \ref{nr:wb-imp}). This minor deviation from machine zero is attributable to standard floating-point round-off accumulation during the algebraic resolution of the linear systems \eqref{eq:fully-discrete-ri-horizontal}-\eqref{eq:fully-discrete-ri-vertical}, rather than any numerical diffusion or truncation error in the scheme. Overall, these results confirm that the exactly well-balanced property of the relaxation framework remains robust and is not compromised by the geometric complexity of the bottom profile.
\begin{table}[htbp]
  \centering
  \caption{Experimental verification of the exactly well-balanced property for the fully explicit schemes (EXP O1 and EXP O2) on a $200 \times 200$ grid. Relative $L^1$-norm errors with respect to the domain size for water height ($h$) and discharges ($q_u, q_v$) at $T = 200$ over the three-mound topography.}
  \label{nr:wb-exp}
  \begin{tabular}{cccccc}
    \toprule
    \multicolumn{3}{c}{EXP O1 PT}  & \multicolumn{3}{c}{EXP O1 TP}                                              \\
    \cmidrule(lr){1-3} \cmidrule(lr){4-6}
    $h$                            & $q_u$                          & $q_v$    & $h$      & $q_u$    & $q_v$    \\
    \midrule
    0.00e+00                       & 0.00e+00                       & 0.00e+00 & 0.00e+00 & 0.00e+00 & 0.00e+00 \\
    \midrule
    \multicolumn{3}{c}{EXP O2 TPT} & \multicolumn{3}{c}{EXP O2 PTP}                                             \\
    \cmidrule(lr){1-3} \cmidrule(lr){4-6}
    $h$                            & $q_u$                          & $q_v$    & $h$      & $q_u$    & $q_v$    \\
    \midrule
    0.00e+00                       & 0.00e+00                       & 0.00e+00 & 0.00e+00 & 0.00e+00 & 0.00e+00 \\
    \bottomrule
  \end{tabular}
\end{table}
\begin{table}[htbp]
  \centering
  \caption{Experimental verification of the exactly well-balanced property for the semi-implicit schemes (IMP O1 and IMP O2) on a $200 \times 200$ grid. Relative $L^1$-norm errors with respect to the domain size at $T = 200$ over the three-mound topography.}
  \label{nr:wb-imp}
  \begin{tabular}{cccccc}
    \toprule
    \multicolumn{3}{c}{IMP O1 PT}  & \multicolumn{3}{c}{IMP O1 TP}                                              \\
    \cmidrule(lr){1-3} \cmidrule(lr){4-6}
    $h$                            & $q_u$                          & $q_v$    & $h$      & $q_u$    & $q_v$    \\
    \midrule
    3.38e-14                       & 1.32e-13                       & 1.31e-13 & 2.96e-14 & 1.45e-13 & 1.43e-13 \\
    \midrule
    \multicolumn{3}{c}{IMP O2 TPT} & \multicolumn{3}{c}{IMP O2 PTP}                                             \\
    \cmidrule(lr){1-3} \cmidrule(lr){4-6}
    $h$                            & $q_u$                          & $q_v$    & $h$      & $q_u$    & $q_v$    \\
    \midrule
    8.50e-15                       & 3.79e-14                       & 3.97e-14 & 6.76e-15 & 8.19e-15 & 8.78e-15 \\
    \bottomrule
  \end{tabular}
\end{table}

\subsection{Order of accuracy test}

To verify the empirical order of accuracy of the proposed schemes, we consider a smooth perturbation of a water-at-rest equilibrium \cite{caballero2024semi}. The initial velocity field is set to zero, $u \equiv v \equiv 0$, and the initial free surface elevation is defined as
\begin{displaymath}
  \eta(x,y,0) = \begin{cases}
    0.05 \! \left(\! 1 + \cos \left(\! 2\pi\dfrac{\sqrt{(x-7000)^2 + (y-7000)^2} - 4750}{3500} \! \right) \! \right) & \text{if }  (x,y) \in C, \\
    0                                                                                                                & \text{otherwise},
  \end{cases}
\end{displaymath}
where $C = \{(x,y) \in \mathbb{R}^2 : |\sqrt{(x-7000)^2 + (y-7000)^2}-4750| < 1750\}$ and the bottom topography is prescribed by
\begin{displaymath}
  H(x,y) = 50 - e^{-\frac{(x-7000)^2}{10^6}}.
\end{displaymath}
The computational domain is $\Omega =  [0,14000]\times[0,14000]$. The numerical errors are evaluated at the final time $T=1$. The stability parameters are strictly set to $\text{CFL}_p = \text{CFL}_t = 1$ for the explicit schemes, whereas for the semi-implicit variants, a significantly relaxed acoustic Courant number of $\text{CFL}_p = 5$ is employed alongside $\text{CFL}_t = 1$.

In the absence of an exact analytical solution, the numerical errors are computed with respect to a high-resolution reference solution obtained on a $3200 \times 3200$ grid. The relative $L^1$-norm errors, normalized by the domain size, are summarized in Tables \ref{tab:accuracy_EXPO1PT} and \ref{tab:accuracy_IMPO1PT} for the first-order schemes, and in Tables \ref{tab:accuracy_EXPO2TPT} and \ref{tab:accuracy_IMPO2TPT} for their second-order counterparts. As expected, both the explicit and semi-implicit formulations successfully achieve their expected order of accuracy.

\begin{table}[htbp]
  \centering
  \caption{Relative $L^1$-norm errors and empirical order of accuracy for the water height ($h$) and discharges ($q_u, q_v$), obtained with the first-order explicit scheme (EXP O1 PT) at time $T=1$.}
  \label{tab:accuracy_EXPO1PT}
  \begin{tabular}{ccccccc}
    \toprule
    \multirow{2}{*}{\begin{tabular}[c]{@{}c@{}}No. of cells\\ ($N_x \times N_y$)\end{tabular}} & \multicolumn{2}{c}{$h$} & \multicolumn{2}{c}{$q_u$} & \multicolumn{2}{c}{$q_v$}                            \\ \cmidrule(lr){2-3} \cmidrule(lr){4-5} \cmidrule(lr){6-7}
                                                                                               & Error                   & Order                     & Error                     & Order & Error    & Order \\ \midrule
    $100 \times 100$                                                                           & 2.24e-04                &                           & 1.15e-04                  &       & 1.14e-04 &       \\
    $200 \times 200$                                                                           & 6.99e-05                & 1.680                     & 4.12e-05                  & 1.478 & 4.09e-05 & 1.479 \\
    $400 \times 400$                                                                           & 2.44e-05                & 1.518                     & 1.58e-05                  & 1.384 & 1.57e-05 & 1.384 \\
    $800 \times 800$                                                                           & 9.13e-06                & 1.422                     & 6.02e-06                  & 1.392 & 5.97e-06 & 1.392 \\
    \bottomrule
  \end{tabular}
\end{table}

\begin{table}[htbp]
  \centering
  \caption{Relative $L^1$-norm errors and empirical order of accuracy for the water height ($h$) and discharges ($q_u, q_v$), obtained with the first-order semi-implicit scheme (IMP O1 PT) at time $T=1$.}
  \label{tab:accuracy_IMPO1PT}
  \begin{tabular}{ccccccc}
    \toprule
    \multirow{2}{*}{\begin{tabular}[c]{@{}c@{}}No. of cells\\ ($N_x \times N_y$)\end{tabular}} & \multicolumn{2}{c}{$h$} & \multicolumn{2}{c}{$q_u$} & \multicolumn{2}{c}{$q_v$}                            \\ \cmidrule(lr){2-3} \cmidrule(lr){4-5} \cmidrule(lr){6-7}
                                                                                               & Error                   & Order                     & Error                     & Order & Error    & Order \\ \midrule
    $100 \times 100$                                                                           & 2.24e-04                &                           & 1.23e-04                  &       & 1.23e-04 &       \\
    $200 \times 200$                                                                           & 6.99e-05                & 1.682                     & 4.57e-05                  & 1.436 & 4.54e-05 & 1.436 \\
    $400 \times 400$                                                                           & 2.43e-05                & 1.521                     & 1.79e-05                  & 1.349 & 1.78e-05 & 1.349 \\
    $800 \times 800$                                                                           & 9.03e-06                & 1.431                     & 6.94e-06                  & 1.367 & 6.89e-06 & 1.367 \\
    \bottomrule
  \end{tabular}
\end{table}

\begin{table}[htbp]
  \centering
  \caption{Relative $L^1$-norm errors and empirical order of accuracy for the water height ($h$) and discharges ($q_u, q_v$), obtained with the second-order explicit scheme (EXP O2 TPT) at time $T=1$.}
  \label{tab:accuracy_EXPO2TPT}
  \begin{tabular}{ccccccc}
    \toprule
    \multirow{2}{*}{\begin{tabular}[c]{@{}c@{}}No. of cells\\ ($N_x \times N_y$)\end{tabular}} & \multicolumn{2}{c}{$h$} & \multicolumn{2}{c}{$q_u$} & \multicolumn{2}{c}{$q_v$}                            \\ \cmidrule(lr){2-3} \cmidrule(lr){4-5} \cmidrule(lr){6-7}
                                                                                               & Error                   & Order                     & Error                     & Order & Error    & Order \\ \midrule
    $100 \times 100$                                                                           & 2.28e-04                &                           & 8.27e-05                  &       & 8.21e-05 &       \\
    $200 \times 200$                                                                           & 5.66e-05                & 2.003                     & 2.20e-05                  & 1.904 & 2.19e-05 & 1.904 \\
    $400 \times 400$                                                                           & 1.40e-05                & 2.015                     & 5.61e-06                  & 1.971 & 5.61e-06 & 1.971 \\
    $800 \times 800$                                                                           & 3.35e-06                & 2.069                     & 1.35e-06                  & 2.064 & 1.34e-06 & 2.064 \\
    \bottomrule
  \end{tabular}
\end{table}

\begin{table}[htbp]
  \centering
  \caption{Relative $L^1$-norm errors and empirical order of accuracy for the water height ($h$) and discharges ($q_u, q_v$), obtained with the second-order semi-implicit scheme (IMP O2 TPT) at time $T=1$.}
  \label{tab:accuracy_IMPO2TPT}
  \begin{tabular}{ccccccc}
    \toprule
    \multirow{2}{*}{\begin{tabular}[c]{@{}c@{}}No. of cells\\ ($N_x \times N_y$)\end{tabular}} & \multicolumn{2}{c}{$h$} & \multicolumn{2}{c}{$q_u$} & \multicolumn{2}{c}{$q_v$}                            \\ \cmidrule(lr){2-3} \cmidrule(lr){4-5} \cmidrule(lr){6-7}
                                                                                               & Error                   & Order                     & Error                     & Order & Error    & Order \\ \midrule
    $100 \times 100$                                                                           & 2.28e-04                &                           & 8.27e-05                  &       & 8.21e-05 &       \\
    $200 \times 200$                                                                           & 5.66e-05                & 2.003                     & 2.21e-05                  & 1.900 & 2.20e-05 & 1.900 \\
    $400 \times 400$                                                                           & 1.40e-05                & 2.015                     & 5.61e-06                  & 1.975 & 5.61e-06 & 1.975 \\
    $800 \times 800$                                                                           & 3.35e-06                & 2.069                     & 1.33e-06                  & 2.084 & 1.32e-06 & 2.084 \\
    \bottomrule
  \end{tabular}
\end{table}

These numerical results demonstrate two key advantages of the proposed framework. First, both the explicit and semi-implicit schemes consistently achieve their orders of accuracy, with the second-order discretizations asymptotically converging at the nominal rate of $2.0$ across all conserved variables ($h, q_u, q_v$). Second, and most notably, the semi-implicit variants yield error magnitudes virtually identical to those of the fully explicit schemes, despite operating at a fivefold larger acoustic Courant number ($\text{CFL}_p = 5$ vs. $\text{CFL}_p = 1$). This confirms that the implicit pressure-splitting strategy successfully bypasses the severe acoustic time-step restriction without introducing excess numerical dissipation or compromising spatio-temporal accuracy.

\subsection{Subcritical regime}

To evaluate the computational efficiency of the proposed schemes in low-Froude-number regimes, we perform a test in a subcritical flow where implicit pressure integration is expected to yield significant performance gains. The test case is adapted from \cite{del2024well} and is posed on the computational domain $\Omega = [0,20] \times [0,20]$. The fluid is initially at rest, $u \equiv v \equiv 0$, with the initial water height given by
\begin{displaymath}
  h(x,y,0) = 2-0.1e^{-((x-10)^2 + (y-10)^2)},
\end{displaymath}
and the bottom is prescribed by
\begin{displaymath}
  H(x,y) = 1.9 + 0.01e^{-((x-10)^2 + (y-10)^2)}.
\end{displaymath}

Simulations are run until a final time of $T=1$. The Froude number is smaller than 0.0012 ($F_r < 1.2 \cdot 10^{-3}$), which makes this example representative of subcritical flows. We fix $\text{CFL}_t=1$ across all configurations, enforce a restrictive $\text{CFL}_p=1$ for the explicit schemes, and significantly relax it to $\text{CFL}_p=8$ for the semi-implicit variants. Table \ref{tab:speed-table} reports the corresponding CPU execution times and the resulting speedup factors on a $1200 \times 1200$ grid.

\begin{table}[htbp]
  \centering
  \caption{Computational efficiency and relative CPU execution time ratios ($R = \text{Time}_{\text{EXP}} / \text{Time}_{\text{IMP}}$) for the subcritical flow benchmark at final time $T = 1.0$ on a uniform $1200 \times 1200$ grid. The explicit configurations, alongside the standard HLL and HLLC-type solvers (previously introduced in Test \ref{test:square-dam-break}), operate under a restrictive acoustic Courant number of $\text{CFL}_p = \text{CFL}_t = 1$. In contrast, the proposed semi-implicit schemes utilize an enlarged pressure Courant number of $\text{CFL}_p = 8$ (while keeping $\text{CFL}_t = 1$). Ratios $R > 1.0$ highlight the performance speedup of the semi-implicit approach over the explicit formulations. Total wall-clock execution times (in seconds) are reported in the rightmost column for explicit schemes and in the bottom row for the semi-implicit methods.}
  \label{tab:speed-table}
  \begin{tabular}{l | c c | c}
    \toprule
    \small{EXP / IMP} & IMP O1 PT & IMP O2 TPT & EXP CPU time (s) \\
    \midrule
    EXP O1 PT         & 5.62      & 2.09       & 325.39           \\
    HLLC-type         & 7.54      & 2.80       & 436.76           \\
    HLL               & 8.46      & 3.15       & 490.24           \\
    EXP O2 TPT        & 16.06     & 5.97       & 930.21           \\
    \midrule
    IMP CPU time (s)  & 57.92     & 155.73     & --               \\
    \bottomrule
  \end{tabular}
\end{table}

The numerical results reported in Table \ref{tab:speed-table} confirm that the
semi-implicit formulations achieve substantial reductions in CPU runtime
compared to both the fully explicit scheme and the standard Riemann
solvers (HLL and HLLC-type), successfully overcoming the acoustic CFL
restriction. When compared at the same order of accuracy, the
semi-implicit scheme is $5.62$ times faster than its explicit counterpart
at first order, and $5.97$ times faster at second order.

The speed-up becomes even more pronounced when the semi-implicit scheme is
compared against the standard HLL and HLLC-type solvers: the first-order
semi-implicit scheme runs between $7.5$ and $8.5$ times faster than these
widely used explicit solvers, reaching a maximum ratio of $8.46$
against the HLL solver. Remarkably, even the second-order semi-implicit
scheme, despite its higher intrinsic cost, remains more than twice as fast
as both standard solvers ($2.80$ against HLLC-type and $3.15$
against HLL).

\subsection{Wave propagation}

Finally, to evaluate the performance of the proposed schemes in a scenario closer to realistic geophysical applications, we consider a wave-propagation benchmark designed to simulate non-linear coastal wave shoaling dynamics, inspired by the numerical studies in \cite{park2007numerical}. The basin is defined on the computational domain $\Omega = [-100,100]\times[0,100]$ with a hyperbolic-tangent bathymetry given by
\begin{displaymath}
  H(x,y)=50 - 45\tanh \left(\frac{-(50 -x)}{20} \right),
\end{displaymath}
and an initial free surface elevation prescribed by
\begin{displaymath}
  \eta(x,y,0) = e^{-\frac{(x-30.25)^2}{10}-\frac{(y-50.25)^2}{20}}.
\end{displaymath}
We set $q_v(x,y,0) \equiv 0$ and $q_u(x,y,0)=100\,\eta(x,y,0)$. To prevent unphysical boundary reflections, open (outflow) boundary conditions are imposed along all domain boundaries. The transport Courant number is fixed to $\text{CFL}_t=1$. For the pressure subsystem, we set $\text{CFL}_p=1$ for explicit computations and a relaxed value of $\text{CFL}_p=6$ for semi-implicit simulations.

\begin{figure}[htbp]
  \centering
  \includegraphics[width=\linewidth]{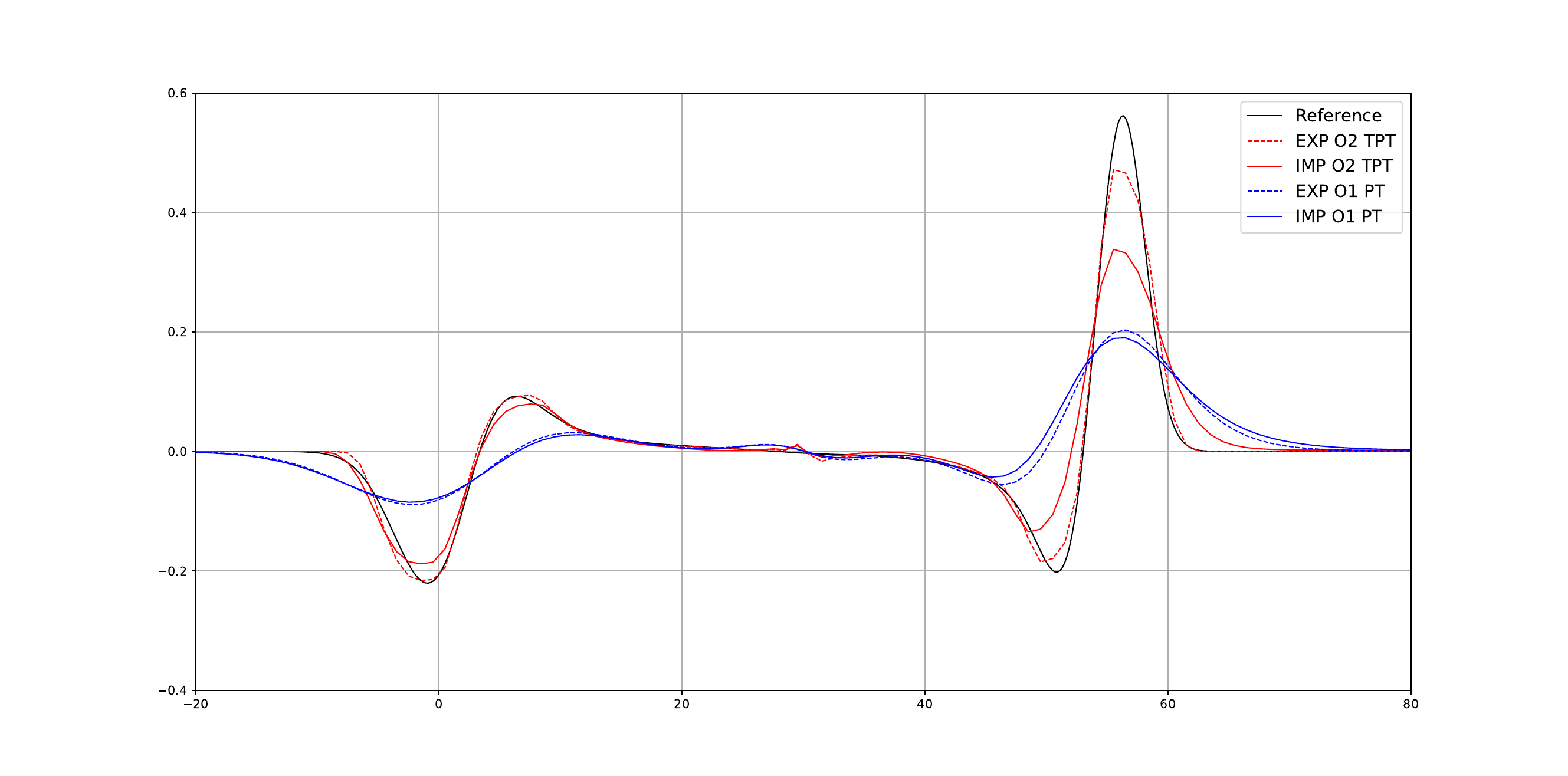}
  \caption{Cross-sectional profile of the free surface elevation $\eta = h - H$ along $y = 50.25$ at time $T = 1$ for the wave propagation benchmark over hyperbolic-tangent bathymetry. Numerical solutions obtained with the proposed first- and second-order explicit ($\text{EXP}$, $\text{CFL}_p = 1$) and semi-implicit ($\text{IMP}$, $\text{CFL}_p = 6$) schemes on a $200\times200$ grid with open boundary conditions are compared against a high-resolution ($1000\times1000$) HLL reference solution. The semi-implicit framework accurately captures non-linear wave shoaling and front propagation without spurious oscillations or artificial boundary reflections.}
  \label{fig:cut_propagation}
\end{figure}

All computations are performed on a $200\times200$ uniform mesh and compared against a high-resolution reference profile produced using an HLL scheme (see Test \ref{test:square-dam-break}) on a $1000\times1000$ grid. Figure \ref{fig:cut_propagation} displays a 1D cross-section of the free surface elevation along $y = 50.25$ at time $T = 1$, providing a direct comparison of the resolving capabilities and wave-front fidelity across the different schemes.

The numerical results indicate that the first-order semi-implicit scheme (IMP O1 PT) clearly outperforms its explicit counterpart by significantly reducing numerical dissipation at the advancing wave front. Although the second-order explicit scheme (EXP O2 TPT) exhibits slightly sharper resolution than the second-order semi-implicit variant (IMP O2 TPT), the proposed semi-implicit framework offers an exceptional accuracy-to-cost trade-off: it delivers an accurate and robust description of complex non-linear wave shoaling, while allowing an increase from $\text{CFL}_p = 1$ to $\text{CFL}_p = 6$ in the permissible time step through relaxed acoustic constraints.

\section{Conclusions}

In this work, we have proposed a novel semi-implicit, exactly well-balanced relaxation scheme for the two-dimensional shallow water equations. This methodology is based on a flux-splitting technique combined with a Suliciu-type relaxation framework, achieving an effective decoupling of acoustic and transport phenomena.

This numerical framework is particularly advantageous for simulating flows in low-Froude-number regimes. In these scenarios, classical fully explicit methods are severely bottlenecked by restrictive time-step limits dictated by fast acoustic waves. By isolating the acoustic component, our approach successfully overcomes this fundamental limitation.

The proposed schemes offer several significant advantages. Primarily, by treating the pressure system implicitly, the method drastically relaxes the standard Courant--Friedrichs--Lewy (CFL) condition. This allows for the use of much larger time steps, yielding a highly efficient computational tool without compromising stability. Furthermore, the scheme exactly preserves water-at-rest steady states by applying a hydrostatic correction exclusively to the pressure system, which is crucial for preventing spurious numerical oscillations in near-equilibrium states over complex bathymetries.

Finally, several extensions of this work are possible. We intend to explore the adaptation of this semi-implicit well-balanced framework to other hyperbolic balance laws, as well as the implementation of higher-order spatial and temporal reconstructions to further improve the scheme's accuracy.

\backmatter

\section*{Declarations}

\bmhead{Funding}
This work has been supported by grant PID2022-137637NB-C21,
funded by MICIU/AEI/10.13039/501100011033 and by ERDF/EU; by grant
PID2025-170651NB-C21, funded by MICIU/AEI/10.13039/501100011033; and by
grant PCI2024-155061-2, funded by MICIU/AEI/10.13039/501100011033 and by
the Deutsche Forschungsgemeinschaft (DFG, German Research Foundation)
through the DFG--AEI 2023 joint call.

\bmhead{Competing interests}
The authors have no competing interests to declare that are relevant to the content of this article.

\bmhead{Data availability}
The data supporting the findings of this study are generated by the numerical
code described in the article and are available from the corresponding author
upon reasonable request.

\bmhead{Author contributions}
All authors contributed equally to the conception and design of the numerical
schemes, the theoretical analysis and the numerical experiments. Jos\'e Ignacio
Ram\'irez-Fuentes was primarily responsible for the implementation of the schemes
and for carrying out the numerical simulations. All authors read and approved the
final manuscript.

\bibliography{references}

\end{document}